\documentclass[a4paper,11pt]{article}

\usepackage{color}
\usepackage{array}
\usepackage{geometry}
\usepackage[english,frenchb]{babel}
\usepackage[T1]{fontenc}
\usepackage{amssymb}
\usepackage{MnSymbol}
\usepackage{mathrsfs}
\usepackage{url}
\usepackage{stmaryrd}
\usepackage[utf8]{inputenc}
\usepackage{graphicx}
\usepackage{pgf,tikz}
\usepackage{subfigure}
\usepackage{caption}
\usepackage{amsmath}
\usepackage{amsfonts}
\usepackage{amsthm}
\usepackage{dsfont}
\usepackage{mathenv}
\usepackage{cancel}
\usepackage{appendix}
\usepackage{multirow}
\usepackage{enumerate}
\usepackage{textcomp}
\usepackage{enumitem}

\usepackage{hyperref}
\usepackage{titletoc}
\titlecontents{chapter}[1.2cm]{}{\bfseries \makebox[0pt][r]{\large\thecontentslabel \enspace }\large \normalfont}{}{\dotfill\contentspage}[]

\usepackage{cleveref}
\newtheorem{theorem}{Theorem}[section]
\newtheorem{prop}[theorem]{Proposition}
\crefname{prop}{proposition}{propositions}
\Crefname{prop}{Proposition}{Propositions}

\crefname{conj}{conjecture}{conjectures}
\Crefname{conj}{Conjecture}{Conjectures}
\newtheorem{remark}[theorem]{Remark}
\crefname{remark}{Remark}{remarks}
\Crefname{remark}{Remark}{Remarks}
\newtheorem{lemma}[theorem]{Lemma}
\newtheorem{cor}[theorem]{Corollary}

\theoremstyle{definition}

\usepackage[utf8]{inputenc}
\usepackage{pgf,tikz}
\usepackage{mathrsfs}
\usetikzlibrary{arrows}
\definecolor{qqzzcc}{rgb}{0.,0.6,0.8}
\definecolor{ffccww}{rgb}{1.,0.8,0.4}
\definecolor{ccqqqq}{rgb}{0.8,0.,0.}
\definecolor{qqzzqq}{rgb}{0.,0.6,0.}

\definecolor{amethyst}{rgb}{0.6, 0.4, 0.8}

\def\Bk{\color{black}}

\newcommand{\R}{\mathbb{R}}

\newcommand{\N}{\mathbb{N}}
\newcommand{\A}{{\mathcal A}}

\newcommand{\x}{\mathbf{x}}

\newcommand{\y}{\mathbf y}
\newcommand{\z}{\mathbf{z}}

\newcommand{\opt}{\mathbf{O}}

\newcommand{\tpp}{, \dots, }

\newcommand{\psibf}{\boldsymbol \psi}

\newcommand{\rb}{\overline r}
\newcommand{\dens}{u}
\newcommand{\densbf}{\boldsymbol \dens}

\newcommand{\dd}{\mu}

\def\epsilon{\varepsilon}
\def\ds{\displaystyle}

\def\epsilon{\varepsilon}
\def\rb{\overline{r}}

\def\tu{\widetilde u}
\def\tN{\widetilde N}

\def\lp {\left( }
\def\rp {\right) }
\newcommand{\be}{\begin{equation}}
\newcommand{\ee}{\end{equation}}
\newcommand{\baa}{\begin{array}}
\newcommand{\eaa}{\end{array}}
\newcommand{\bi}{\begin{itemize}}
\newcommand{\ei}{\end{itemize}}
\newcommand{\rmax}{r_{\max}}

\newcommand{\md}{\mathrm{d}}

\title{\bf{Adaptation in a heterogeneous environment. \\I: Persistence versus extinction}\thanks{The project leading to this publication has received funding from Excellence Initiative of Aix-Marseille University~-~A*MIDEX, a French ``Investissements d'Avenir'' programme, and from the ANR project RESISTE (ANR-18-CE45-0019).}}

\author{F. Hamel$^{\hbox{\small{ a}}}$, F. Lavigne$^{\hbox{\small{ a,b,c}}}$, and L. Roques$^{\hbox{\small{ b}}}$,   \\
\\
\footnotesize{$^{\hbox{a }}$Aix Marseille Univ, CNRS, Centrale Marseille, I2M, Marseille, France}\\
\footnotesize{$^{\hbox{b }}$INRAE, BioSP, 84914, Avignon, France}\\
\footnotesize{$^{\hbox{c }}$LAMAV, Universit\'e Polytechnique des Hauts-de-France, France}
}

\date{}

\begin{document}

\selectlanguage{english}

\maketitle

\begin{abstract}
Understanding how a diversity of plants in agroecosystems affects the adaptation of pathogens is a key issue in agroecology. We analyze PDE systems describing the dynamics of adaptation of two phenotypically structured populations, under the effects of mutation, selection and migration in a two-patches environment, each patch being associated with a different phenotypic optimum. We consider two types of growth functions that depend on the $n-$dimensional phenotypic trait: either local and linear or nonlocal nonlinear. In both cases, we obtain existence and uniqueness results as well as a characterization of the large-time behaviour of the solution (persistence or extinction) based on the sign of a principal eigenvalue. We show that migration between the two environments decreases the chances of persistence, with in some cases a ‘lethal migration threshold’ above which persistence is not possible. Comparison with stochastic individual-based simulations shows that the PDE approach accurately captures this threshold. Our results illustrate the importance of cultivar mixtures for disease prevention and control.
\end{abstract}

\section{Introduction}

Phenotypic differences between populations generally appear as a consequence of differential selection regimes \cite{Orr98}. For instance, in the absence of migration, the adaptation of a population to local habitat conditions leads to a particular phenotypic distribution. In asexual populations, a standard way to describe the gene -- environment interaction is to use Fisher's geometrical model (FGM) \cite{MarLen15,Ten14}. In this approach, each individual in the population is characterized by a multivariate phenotype at a set of $n$ traits, i.e., a vector $\x \in \R^n$. This vector $\x$ determines the fitness $r(\x)$ (the reproductive success of the individual) through its quadratic distance with respect to an optimum $\opt \in \R^n$ associated with the considered environment:
$$r(\x)=\rmax-\frac{\|\x-\opt\|^2}{2},$$
with $\rmax>0$ the fitness of the optimal phenotype. Throughout the paper, $\|\cdot \|$ denotes the Euclidean norm in $\R^n$.

\paragraph{PDE models.} Under the assumption of the FGM, recent models of asexual adaptation based on partial differential equations (PDE) \cite{AlfCar17,AlfVer18,HamLavMarRoq20} typically describe the dynamics of the phenotype distribution $q$ of a population in a single environment, with equations of the form:
$$\forall t >0, \, \forall \x\in \R^n,\quad \partial_t q(t,\x)=\frac{\mu^2}{2}\Delta q(t,\x)+[r(\x)-\rb(t)]\, q(t,\x),$$where the Laplace operator describes the mutation effects on the phenotype (see~\cite[Appendix]{HamLavMarRoq20} for the derivation of this term in this framework), and the term $[r(\x)-\rb(t)]\, q(t,\x)$ describes the effects of selection~\cite{TsiLev96}, with $\rb(t)=\int_{\R^n}r(\x)\,q(t,\x)\,\md\x$ the mean fitness in the population at time $t$. Extensions to temporally changing environments (with an optimum $\opt(t)$) have also been proposed \cite{RoqPat20}. In all those cases, it was possible to describe the full dynamics of adaptation, by deriving explicit expressions for~$\rb(t)$. 

The simplest case with a unique constant optimum $\opt$ describes the adaptation of a population to an abrupt environmental change,  e.g., when a bacterial population is exposed to an antibiotic or during a host shift for a virus. The theory in \cite{HamLavMarRoq20,MarRoq16} shows that the mean growth rate satisfies $\rb(t)=\rmax-\mu\, n/2 \, \tanh(\mu \, t)+ (\rb(0)-\rmax)/\cosh(\mu \, t)^2$ if one starts from a clonal population with initial growth rate $\rb(0)$. Thus, $\rb(t)\to \rmax-\mu\, n/2$ at large times, independently of the initial phenotype distribution. The quantity $\mu\, n/2$ corresponds to the amount of maladaptation due to mutations, or in other terms the \emph{mutation load}. If the mutation load exceeds the fitness of the optimal phenotype $\rmax$, the population is doomed to extinction, as $\rb(t)<0$ at large times.  Lethal mutagenesis theory is based on this expected effect of high mutation rates on the viability of pathogens (and especially viruses) \cite{AncLam19,BulSan07,BulWil08}.

Agroecological theory rather relies on the diversification of host species to reduce the viability of pathogens \cite{CaqGas20}. Here, we consider a spatially heterogeneous environment, made of two habitats, each of them corresponding to a different phenotype optimum, $\opt_1$ and $\opt_2$. The main issue that we are going to deal with is to determine the respective effect of the migration between the two habitats and of the phenotypic distance between the two habitats on the fate (persistence or extinction) of the total population. The underlying question is to determine whether a system with two hosts, corresponding to two species or two genetic variants of the same species, connected by migration is more resilient to invasion by a pathogen than a single host, and how this depends on the parameter values. This type of question has already been considered with comparable models in \cite{MirGan20}, in a particular regime of parameters such that the effect of the mutation is low, while the mutation rate is large enough, and in dimension $n=1$. The authors have used a specific method based on constrained Hamilton-Jacobi equations (e.g., \cite{BarMir09,DieJab05,GanMir17,LorMir11,PerBar08} for more details on this method), to find an accurate analytic approximation of the equilibrium phenotype distribution and the population size in each habitat. They found that, when the two environments are symmetric (same mutation parameters, same selection pressure, same competition intensity and same migration rates), there exists an explicit threshold for the migration rate, which depends on the phenotypic distance between the two habitats. When the migration rate is above this threshold, the two subpopulations are well-mixed so that the total equilibrium population is monomorphic or `generalist'. On the contrary, when the migration is below the threshold, the two subpopulations stay different, causing dimorphism in the phenotype density for the global population: the equilibrium population is made of two `specialists'. They also obtained some results  in the general case, without the symmetry assumption. 

As we focus here on persistence/extinction issues, instead of dealing with the phenotype distribution $q(t,\x)$, we are interested in the \emph{phenotype density} $\dens(t,\x)$, i.e.,  $\dens(t,\x)=q(t,\x) \, N(t)$, with $N(t)$ the population size at time $t$. We therefore deal with systems of the form:
\begin{equation}
\forall\,  t\ge0, \ \forall\,\x\in \R^n,\, \left\{\begin{array}{rcl}
\partial_t\dens_1 (t,\x) & = &\ds \frac{\dd^2} 2 \ \Delta\dens_1(t,\x)+f_1(\x,[\dens_1])+\delta\,[\dens_2(t,\x) - \dens_1(t,\x)], \vspace{2mm}\\
     \partial_t\dens_2 (t,\x) & =&\ds \frac{\dd^2} 2 \ \Delta\dens_2(t,\x)+f_2(\x,[\dens_2]) +\delta\,[\dens_1(t,\x) - \dens_2(t,\x)],
\end{array}\right.  
\label{eq:sys_main1}
\end{equation}
with $\dens_i$ the phenotype density in the habitat $i\in\{1,2\}$, $\delta>0$ the migration rate, and $\mu>0$ a mutational parameter. Note that the migration and mutation parameters are assumed to be identical in the two habitats. This is a simplifying assumption which leads to symmetry properties of the solutions that are important in most of our results. However, some results also deal with the case of non-symmetric systems of the type:
\begin{equation}
\forall\,  t\ge0, \ \forall\,\x\in \R^n,\, \left\{\begin{array}{rcl}
\partial_t\dens_1 (t,\x) & = &\ds \frac{\dd^2} 2 \ \Delta\dens_1(t,\x)+f_1(\x,[\dens_1])-\delta_{1,1}\dens_1(t,\x)+\delta_{1,2}\dens_2(t,\x), \vspace{2mm}\\
\partial_t\dens_2 (t,\x) & =&\ds \frac{\dd^2} 2 \ \Delta\dens_2(t,\x)+f_2(\x,[\dens_2])+\delta_{2,1}\dens_1(t,\x)-\delta_{2,2}\dens_2(t,\x),
\end{array}\right.  
\label{eq:sys_main1bis}
\end{equation}
where the $\delta_{i,j}$'s are positive constants. We however still keep in~\eqref{eq:sys_main1bis} the same mutational parameter~$\dd$ in both habitats: this assumption is used in the derivation of explicit bounds in the resolution of the Cauchy problem (see the proof of Theorem~\ref{th:exist_uniq_dens} in Section~\ref{sec:31} below for more details).

We may assume two different types of growth functions $f_i$. We first state that, in both cases, the fitness of a phenotype $\x$ in the habitat $i$ is given by:
\be\label{def:ri}
r_i(\x)=\rmax-\frac{\|\x-\opt_i\|^2}{2},
\ee
or more generally:
\be\label{def:ribis}
r_i(\x)=r_{\max,i}-\frac{\|\x-\opt_i\|^2}{2},
\ee
for some real numbers $r_{\max,i}$. Notice in particular that the fitnesses $r_i$ are unbounded in the phenotypic space $\R^n$ and, since they are involved in the definition of the growth functions $f_i$ for both types listed below, the systems~\eqref{eq:sys_main1} and~\eqref{eq:sys_main1bis} of unknowns~$(\dens_1,\dens_2)$ then have unbounded coefficients.

\paragraph{The first type $($Malthusian$)$:}
\begin{equation}\label{eq:f_type1}
    f_i(\x,[\dens_i])=r_i(\x) \, \dens_i(t,\x),
\end{equation}
corresponds to the standard assumption of Malthusian population growth, namely:
\begin{equation} \label{eq:ODEmalthus}
\forall \, t\ge 0,\ \     \left\{
\begin{array}{c}
N_1'(t)=\rb_1(t)\, N_1(t) +\delta\,[N_2(t)-N_1(t)],\vspace{3pt}\\
N_2'(t)=\rb_2(t)\, N_2(t) +\delta\,[N_1(t)-N_2(t)],\\
\end{array}
\right.
\end{equation}
in the symmetric case~\eqref{eq:sys_main1}, and to:
\begin{equation} \label{eq:ODEmalthusbis}
\forall \, t\ge 0,\ \     \left\{
\begin{array}{c}
N_1'(t)=\rb_1(t)\, N_1(t)-\delta_{1,1}N_1(t)+\delta_{1,2}N_2(t),\vspace{3pt}\\
N_2'(t)=\rb_2(t)\, N_2(t) +\delta_{2,1}N_1(t)-\delta_{2,2}N_2(t),\\
\end{array}
\right.
\end{equation}
in the general case~\eqref{eq:sys_main1bis}, where in both cases $N_i(t)$ is the population size in habitat $i$ at time $t$:
\be
N_i(t) = \int_{\R^n} \dens_i(t,\x)\, \md\x,
\label{dfn:pop_size_pde}
\ee
and $\rb_i(t)$ the mean fitness of the individuals located in habitat $i$ at time $t$:
\begin{equation}\label{dfn:meanfit}
\rb_i(t)=\frac{1}{N_i(t)}\int_{\R^n}r_i(\x) \,  \dens_i (t,\x) \, \md\x.
\end{equation}
Note that, with $f_i$ of the type~\eqref{eq:f_type1}, the systems~\eqref{eq:sys_main1} and~\eqref{eq:sys_main1bis} are local cooperative systems since the right-hand side of the equation of each component is nondecreasing with respect to the other component, and since the right-hand side only depends on the densities for the phenotype~$\x$. As a consequence, the maximum principle holds for~\eqref{eq:sys_main1} and~\eqref{eq:sys_main1bis} in this first type~\eqref{eq:f_type1}, that is, if $\mathbf{u}=(u_1,u_2)$ and $\mathbf{v}=(v_1,v_2)$ are two classical solutions of~\eqref{eq:sys_main1} or~\eqref{eq:sys_main1bis} which are locally bounded in time and are such that $\mathbf{u}(0,\cdot)\le\mathbf{v}(0,\cdot)$ in $\R^n$ (in the sense of componentwise inequalities), then $\mathbf{u}(t,\cdot)\le\mathbf{v}(t,\cdot)$ in~$\R^n$ for all $t>0$.

\paragraph{The second type $($density-dependent$)$:}
\begin{equation}\label{eq:f_type2}
    f_i(\x,[\dens_i])=\lp r_i(\x) - \int_{\R^n} \dens_i(t,\y) \, \md\y\rp\,\dens_i(t,\x),
\end{equation}
for which the existence and uniqueness of solutions will only be proven for the symmetric system~\eqref{eq:sys_main1} with fitnesses given by~\eqref{def:ri} and with symmetric initial conditions (see especially Theorem~\ref{th:exist_uniq_dens} in Section~\ref{sec:results}), corresponds to the standard assumption of logistic population growth:
\begin{equation} \label{eq:ODElogistic}
\forall\,t\ge 0,\ \      \left\{
\begin{array}{c}
N_1'(t)=\rb_1(t)\, N_1(t) - N_1(t)^2 +\delta\,[N_2(t)-N_1(t)],\vspace{3pt}\\
N_2'(t)=\rb_2(t)\, N_2(t) - N_2(t)^2 +\delta\,[N_1(t)-N_2(t)],
\end{array}
\right.
\end{equation}
with $N_i(t)$ and $\rb_i(t)$ as in~\eqref{dfn:pop_size_pde}-\eqref{dfn:meanfit}. Note that, with $f_i$ of the type~\eqref{eq:f_type2}, the system~\eqref{eq:sys_main1} is a nonlocal system, the nonlocality in the form of an internal competition. As a consequence, the maximum principle does not hold for~\eqref{eq:sys_main1} in this second type~\eqref{eq:f_type2}.

\vspace{0.5cm}

\noindent{}In this work, after the analysis of the Cauchy problem for~\eqref{eq:sys_main1} and~\eqref{eq:sys_main1bis}, we will mainly focus our attention, for problem~\eqref{eq:sys_main1} with fitnesses~\eqref{def:ri}, and for both types of growth functions~\eqref{eq:f_type1} and~\eqref{eq:f_type2}, on the effects of the migration parameter $\delta$ and of the \emph{habitat difference} defined by:
\be\label{def:mD}
m_D := \frac {\|\opt_1-\opt_2\|^2} {2}>0.
\ee

\paragraph{Stochastic model.} Before stating our main results, and in order to underline the interest of the PDE approach, we compare its accuracy with a standard Wright-Fisher individual-based stochastic model (IBM), with mutation, selection and migration.

In this IBM, each individual  is characterized by a phenotype $\x\in\R^n$, and a corresponding fitness $r_i(\x)$, depending on the position (i.e., the habitat $i=1$ or $i=2$) of the individual. The populations in the two habitats are initially clonal (with all of the phenotypes set at $(\opt_1+\opt_2)/2$), and of size $N_i(0)=N^0$. Then, at each time step (the model is discrete in time), the \emph{reproduction-selection step} is simulated by drawing a Poisson number of offspring, for each individual, with rate $\exp(r_i(\x))$ (Darwinian fitness, the discrete-time counterpart of $r_i(\x)$). Then,  the \emph{mutation step} is simulated by randomly drawing, for each individual, a Poisson number of mutations, with rate~$U>0$. Each single mutation has a random phenotypic effect $\md\x\in \R^n$ drawn into a multivariate Gaussian distribution: $\md \x\sim \mathcal{N}(0,\lambda I_n)$, where $\lambda>0$ is the mutational variance at each trait, and~$I_n$ is the identity matrix of size $n\times n$. Multiple mutations in a single individual have additive effects on phenotype. Lastly, the \emph{migration step} consists in sending individuals from the first habitat into the second (resp. from the second into the first): the numbers of migrants are drawn in a Poisson law with parameter $\delta N_1(t)$ (resp. $\delta N_2(t)$), and the migrants are randomly sampled in the populations. 

\begin{figure}
     \centering
     \begin{subfigure}[]
        {\includegraphics[scale=0.48]{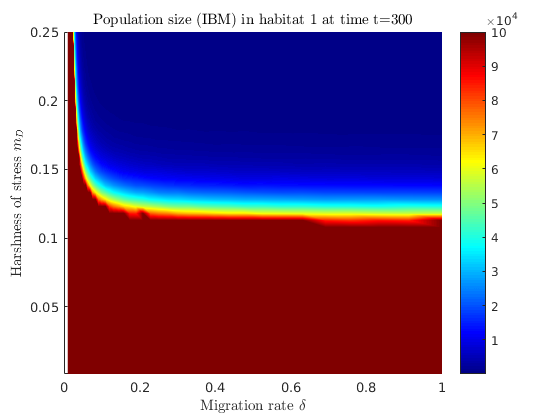}}
     \end{subfigure}
     \begin{subfigure}[]
        {\includegraphics[scale=0.48]{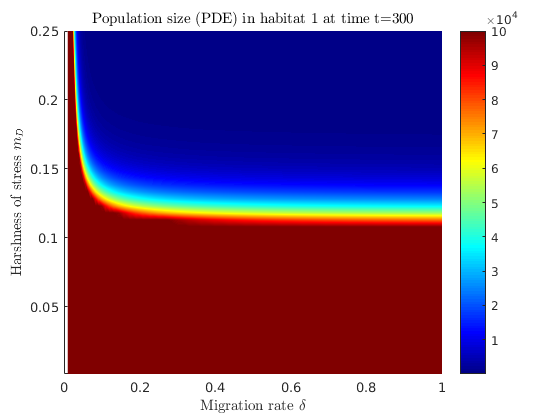}}
     \end{subfigure}
        \caption[Persistence vs extinction: effect of the migration rate $\delta$ and the habitat difference $m_D$]{\textbf{Persistence vs extinction: effect of the migration rate $\delta$ and the habitat difference $m_D$.} Total population size $N_1+N_2$, given (a) by simulation of the stochastic model (average result over 50 replicate simulations) (b) by numerically solving \eqref{eq:sys_main1} with $f_1,f_2$ given by~\eqref{eq:f_type1}. The parameters are $U=1/6$, $\lambda=1/300$, $n=2$, $\rmax=1/18$ and $\mu^2=\lambda \, U$, and the results are computed at $t=300$. Initially, each habitat $i\in\{1,2\}$ has $N^0=10^4$ individuals, all of them with the phenotype $(\opt_1+\opt_2)/2$. }
        \label{fig:graph_mDcrit}
\end{figure}

\paragraph{Numerical comparison between the PDE and stochastic models.}
We simulated the IBM until a time $t=300$, and compared the result with the numerical solution of the PDE model~\eqref{eq:sys_main1}  with $\mu^2=\lambda\, U$  (see~\cite[Appendix]{HamLavMarRoq20} for a justification of this parameter choice), and with the first type of growth function (Malthusian), as the IBM does not take density-dependence into account. The solution of the PDE was computed using the method of lines coupled with the Runge-Kutta ODE solver Matlab\textsuperscript{\tiny\textregistered}  \emph{ode45}.  The results are presented in Fig.~\ref{fig:graph_mDcrit}. We observe a very good agreement between the results obtained with the IBM and the PDE, with in both cases a strong dependence of the persistence/extinction behaviour with respect to the parameters~$\delta$ and~$m_D$. 

\paragraph{Aim of this paper.} A rich theoretical literature has considered the effects of migration on the evolution of sexually reproducing populations in a heterogeneous landscape (e.g., \cite{BolNos07,GarKir97,KirBar97}). Migrants hybridize/recombine with locally adapted genotypes, leading to a locally reduced fitness. This reduction in fitness is known as a \emph{migration load} \cite{GarKir97}. We focus here on asexuals orga\-nisms, such as viruses or bacteria. In this context, although the evolutionary effect of migration is limited compared to the case of sexual organisms, theoretical studies have shown that it also leads to a decrease in local adaptation. These studies analyzed the joint effects of selection and migration (and sometimes mutation) on the adaptation of a population in two-patch models (e.g., \cite{DebRon13,MesCzi97,MirGan20}) or in complex networks of patches interconnected via dispersal (e.g., \cite{PapDav13}).  The general conclusion is that high migration rates favour generalist strategies (i.e., intermediate phenotypes between the optimal ones), and decrease the potential for global persistence over the whole environment, which may have important applications in the management of pathogens. While these results are based on formal analytical computations and numerical simulations, and deal with a unique trait, we intend here to confirm whether they can be derived rigorously, in the above-described $n$-dimensional PDE framework, based on elliptic and parabolic theory. More precisely, our main goal is to set on a firm mathematical basis the behaviour observed in Fig.~\ref{fig:graph_mDcrit}, based on the sign of the principal eigenvalue of a system of linear elliptic equations, and to study the dependence of this eigenvalue with respect to the model parameters. The main results are presented in the next section, and discussed in Section~\ref{sec:discu}.

\section{Main results \label{sec:results}}

Without loss of generality, we assume that the optima $\opt_1$ and $\opt_2$ are located along the $x_1$-axis and are symmetric with respect to the origin, i.e., there exists $\beta\ge0$ such that:
\be \label{def:opt_beta}
\opt_1 = (-\beta,0\tpp 0), \quad \hbox{ and } \quad \opt_2 = (\beta,0\tpp 0).
\ee
This means that the habitat difference $m_D$ defined in~\eqref{def:mD} is equal to:
$$m_D=2\beta^2.$$
We are especially interested in the case $\beta>0$ and $m_D>0$, but the case $\beta=m_D=0$ can also be considered from a mathematical point of view. We also assume in some statements that the two densities $\dens_1$ and $\dens_2$ are initially symmetric with respect to the hyperplane $\big\{\x = (x_1\tpp x_n) \in \R^n, \, x_1=0\big\}$:
\[
\tag{SH}
    \forall\,\x \in \R ^n , \ \ \dens_1^0(\x) = \dens_2^0(\iota(\x))=:u^0(\x),\ \densbf^0(\x):=(\dens_1^0(\x),\dens_2^0(\x))=(u^0(\x),u^0(\iota(\x))),
\label{hyp:HS}
\]
with:
\be
\forall\,\x =(x_1\tpp x_n) \in \R^n,\ \ \iota(\x) = (-x_1,x_2\tpp x_n).
\label{dfn:iota}
\ee

\paragraph{The Cauchy problem.} We first show that the Cauchy problems associated with~\eqref{eq:sys_main1} or~\eqref{eq:sys_main1bis} admit a unique solution, under some assumptions on the initial condition $\densbf^0=(\dens^0_1,\dens^0_2)$:
\begin{enumerate}
    \item[(H1)] $u^0_1,\, u^0_2\in C^{2,\alpha}(\R^n)$ for some $\alpha \in (0,1)$;
    \item[(H2)] $0<\displaystyle N^0_i:=\int_{\R^n}u^0_i(\x)\, \md\x<+\infty$ for $i=1,2$;
    \item[(H3)] there exists a nonincreasing function $g : \R_ + \to \R$ (with $\R_+=[0,+\infty)$) such that:
    \begin{itemize}
        \item[(i)] $0\le u^0_i(\x) \le g(\|\x-\opt_i\|)$ for all $\x \in \R^n$ and $i=1,2$; 
        \item[(ii)] the function $r\mapsto r^{n+1} g(r)$ belongs to $L^1(\R_+)$ and converges to $0$ as $r\to+\infty$.
    \end{itemize}
\end{enumerate}
Hereafter, we always make the assumptions~(H1)-(H3), while~\eqref{hyp:HS} will be assumed only in some statements.

Our first main result provides the existence and uniqueness of the density of phenotypes, for both types~\eqref{eq:f_type1} and~\eqref{eq:f_type2} of growth functions $f_1,f_2$. 

\begin{theorem} \label{th:exist_uniq_dens}
\begin{itemize}
\item[{\rm{(i)}}] Assume that $f_1,f_2$ are both of the first type~\eqref{eq:f_type1}, and that the initial condition $\densbf^0=(u^0_1,u^0_2)$ satisfies~{\rm{(H1)-(H3)}}. Then, there exists a unique solution $\densbf=(\dens_1,\dens_2) \in C^{1,2}(\R_+ \times \R^n) $ of the system~\eqref{eq:sys_main1bis} with fitnesses~\eqref{def:ribis}, such that $\densbf\in L^\infty((0,T)\times\R^n)$ for all $T>0$, $\dens_i(t,\x)\to0$ as $\|\x\|\to+\infty$ locally uniformly in $t\in\R_+$, $\dens_i>0$ in $(0,+\infty)\times\R^n$, and the population sizes $N_i:\R_+\to(0,+\infty)$ are of class $C^1$ and satisfy~\eqref{eq:ODEmalthusbis}, with continuous mean fitnesses $\rb_i:\R_+\to\R$.
\item[{\rm{(ii)}}] Assume that $f_1,f_2$ are both of the first type~\eqref{eq:f_type1}, or both of the second type~\eqref{eq:f_type2}, and that the initial condition $\densbf^0=(u^0,u^0\circ\iota)$ satisfies~\eqref{hyp:HS} and {\rm{(H1)-(H3)}}. Then, there exists a unique solution $\densbf=(\dens_1,\dens_2) \in C^{1,2}(\R_+ \times \R^n) $ of the system~\eqref{eq:sys_main1} with fitnesses~\eqref{def:ri}, such that $\densbf\in L^\infty((0,T)\times\R^n)$ for all $T>0$, $\dens_i(t,\x)\to0$ as $\|\x\|\to+\infty$ locally uniformly in $t\in\R_+$, $\dens_i>0$ in $(0,+\infty)\times\R^n$, the population sizes $N_i:\R_+\to(0,+\infty)$ are of class $C^1$, the mean fitnesses $\rb_i:\R_+\to\R$ are continuous, $\densbf$ is symmetric in the sense that:
\be
\dens_1(t,\x) = \dens_2(t,\iota(\x)),\ \hbox{ for all }t\ge 0 \hbox{ and } \x=(x_1\tpp x_n) \in \R^n,
\label{rel:dens1_dens2:symm}
\ee
with $\iota$ defined in \eqref{dfn:iota}, and:
\be\label{eq:Nrb}
N_1(t) = N_2(t)=:N(t)\ \hbox{ and }\ \rb_1(t) = \rb_2(t)=:\rb(t),\  \hbox{ for all }t\ge0.
\ee
Moreover, the population sizes $N_1=N_2$ satisfy \eqref{eq:ODEmalthus} if $f_1,f_2$ are of the first type~\eqref{eq:f_type1}, whereas $N_1=N_2$ satisfy~\eqref{eq:ODElogistic} if $f_1,f_2$ are of the second type~\eqref{eq:f_type2}. In both cases, the functions~$\dens_i$ satisfy the nonlocal parabolic equation:
\begin{equation}
\partial_t \dens_i(t,\x) = \frac{\dd^2}{2}\, \Delta \dens_i(t,\x) +f_i(\x,[\dens_i])+\delta\,[\dens_i(t,\iota(\x))-\dens_i(t,\x) ],
\label{pde:denssym}
\end{equation}
for all $t\ge0$ and  $\x\in\R^n.$
\end{itemize}
\label{thm:existunidens}
\end{theorem}

\begin{remark}\label{rem:nonsym}{\rm For problem~\eqref{eq:sys_main1} with fitnesses~\eqref{def:ri} and growth functions $f_i$ of the first type~\eqref{eq:f_type1}, if the initial condition $\densbf^0=(u^0_1,u^0_2)$ satisfies~{\rm{(H1)-(H3)}} but not~\eqref{hyp:HS}, then the population sizes~$N_1$ and $N_2$ still satisfy~\eqref{eq:ODEmalthus} but the $N_i$ and $\rb_i$ do not satisfy~\eqref{eq:Nrb} in general and the subsequent analysis then becomes more involved. For problem~\eqref{eq:sys_main1} with fitnesses~\eqref{def:ri} and growth functions $f_i$ of the second type~\eqref{eq:f_type2}, the existence and uniqueness is established only for symmetric solutions satisfying~\eqref{rel:dens1_dens2:symm}, since the proof, which is based on a change of functions amounting to a system with the first type~\eqref{eq:f_type1}, uses a key ingredient, that is the equality of the corresponding population sizes.}
\end{remark}

\paragraph{Persistence vs extinction.} Before going further on, we give a precise meaning to the notions of persistence and extinction. By extinction, we mean that the total population size $N_1(t)+N_2(t)$ satisfies:
$$N_1(t)+N_2(t)\to0\ \hbox{ as $t\to +\infty$}.$$
By persistence, we mean that the population does not get extinct at large times, that is:
$$\limsup_{t\to+\infty}\ \big(N_1(t)+N_2(t)\big)>0.$$
To analyze the effect of the parameter values on the persistence/extinction behaviour of the systems~\eqref{eq:sys_main1} or~\eqref{eq:sys_main1bis}, we consider an eigenvalue problem (see \cite{CanCos03} for several other examples of persistence/extinction results \textit{via} eigenvalue problems in bounded domains).

For any $R>0$, we denote by $\A$ the self-adjoint differential operator:
\begin{equation}\label{eq:defA}
\A:=-\, {\dd^2\over 2} \Delta - \left( \begin{matrix}
r_1(\x) -\delta & \delta \\
\delta & r_2(\x) - \delta
\end{matrix} \right),
\end{equation}
for problem~\eqref{eq:sys_main1} with fitnesses~\eqref{def:ri}, and:
\begin{equation}\label{eq:defAbis}
\A:=-\, {\dd^2\over 2} \Delta - \left( \begin{matrix}
r_1(\x) -\delta_{1,1} & \delta_{1,2} \\
\delta_{2,1} & r_2(\x) - \delta_{2,2}
\end{matrix} \right),
\end{equation}
for problem~\eqref{eq:sys_main1bis} with fitnesses~\eqref{def:ribis}, acting here, in both cases, on functions in $[W^{2,n}_{loc}(B(0,R))\cap C_0(\overline {B(0,R)})]^2$, with $B(0,R)$ the open Euclidean ball of~$\R^n$ of center $0$ and radius $R>0$, and $C_0(\overline{B(0,R)})$ the space of continuous functions in~$\overline{B(0,R)}$ which vanish on $\partial B(0,R)$. It follows from \cite[Theorem~13.1]{BusSir04} or \cite[Theorem 1.1]{Sweers92} that the operator~$\A$ defined in~\eqref{eq:defA}, resp.~\eqref{eq:defAbis}, admits a unique principal eigenvalue $\lambda^{R}\ge -\rmax$, resp.:
\be\label{ineqlambda}
\lambda^R\ge\min(-r_{\max,1}+\delta_{1,1}-\delta_{1,2},-r_{\max,2}+\delta_{2,2}-\delta_{2,1}),
\ee
and a unique (up to multiplication by a positive constant) pair of positive (in $B(0,R)$) eigenfunctions $(\varphi_1^{R}, \varphi_2^{R})\in [W^{2,n}_{loc}(B(0,R))\cap C_0(\overline{B(0,R)})]^2$, satisfying:
$$\A\, (\varphi_1^{R},\varphi_2^{R})=\lambda^{R}\, (\varphi_1^{R}, \varphi_2^{R})\ \hbox{ in }B(0,R).$$
Moreover, the functions $\varphi_i^{R}$ are of class $C^\infty_0(\overline{B(0,R)})=C^\infty(\overline{B(0,R)})\cap C_0(\overline{B(0,R)})$ by standard elliptic estimates, and the eigenvalue~$\lambda^{R}$ is characte\-rized by the following minmax formula:
\[
\lambda^{R} = \sup\limits_{(\psi_1,\psi_2) \in E}\quad \inf\limits_{\x \in B(0,R),\,i\in\{1,2\}} \quad  {(\A\, (\psi_1,\psi_2))_i(\x) \over \psi_i(\x)},
\]
with:
\[
E=\big\{(\psi_1,\psi_2)\in [C^2(B(0,R))\cap C(\overline{B(0,R)})]^2,\ \psi_i(\x)>0\hbox{ for all }\x\in B(0,R)\hbox{ and }i\in\{1,2\}\big\}.
\]
This formula readily implies that the map $R\mapsto \lambda^{R}$ is nonincreasing. Since $\lambda^{R}$ is bounded from below independently of $R\ge1$, the quantity $\lambda^{R}$ admits a finite limit as $R\to +\infty$:
\begin{equation}\label{eq:def_lambda}
\lambda:=\lim_{R\to +\infty} \lambda^{R},
\end{equation}
and:
\be\label{ineqlambdabis}
\left\{\baa{l}
\lambda\ge-\rmax\hbox{ for~\eqref{eq:defA}},\vspace{3pt}\\
\lambda\ge\min(-r_{\max,1}+\delta_{1,1}-\delta_{1,2},-r_{\max,2}+\delta_{2,2}-\delta_{2,1})\hbox{ for~\eqref{eq:defAbis}}.\eaa\right.
\ee

\begin{remark}\label{rem:phi_symm}{\rm For the operator~$\A$ in~\eqref{eq:defA} with fitnesses~\eqref{def:ri}, the eigenfunctions $\varphi_1^{R}$ and $\varphi_2^{R}$ also satisfy a symmetry property, namely $\varphi_1^{R}(\x)=\varphi_2^{R}(\iota(\x))$ for all $\x\in\overline{B(0,R)}$, with $\iota$ defined in \eqref{dfn:iota}. Indeed, by setting $\widetilde \varphi_1 (\x)=\varphi_2^{R}(\iota(\x))$ and $\widetilde \varphi_2 (\x)=\varphi_1^{R}(\iota(\x))$ for $\x\in\overline{B(0,R)}$, and by using the symmetry assumption~\eqref{def:opt_beta}, one has $\A\,(\widetilde \varphi_1, \widetilde \varphi_2)=\lambda^{R}\, (\widetilde \varphi_1, \widetilde \varphi_2)$. By uniqueness (up to multiplication) of the pair of principal eigenfunctions, there exists $K>0$ such that $(\widetilde \varphi_1, \widetilde \varphi_2)=K \, (\varphi_1^{R}, \varphi_2^{R})$. At $\x=0$, we get that $\varphi_2^{R}(0)=\widetilde \varphi_1(0)=K\,\varphi_1^{R}(0)$ and $\varphi_1^{R}(0)=\widetilde \varphi_2(0)=K\,\varphi_2^{R}(0)$. Therefore $K=1$ and $\varphi_1^{R}=\varphi_2^{R}\circ\iota$ in $\overline{B(0,R)}$. Thus, for the operator $\A$ in~\eqref{eq:defA} with fitnesses~\eqref{def:ri}, we may (and we have to) take the same normalization
condition for $\varphi_1^{R}$ and $\varphi_2^{R}$. In the proofs, we either assume that $\|\varphi_i^{R}\|_{L^1(B(0,R))} =1$ for both $i\in\{1,2\}$, or $\varphi_i^{R}(0) =1$ for both $i\in\{1,2\}$.}
\end{remark}

The large time behaviour of the population size is closely related to the sign of the quantity~$\lambda$ defined in~\eqref{eq:def_lambda}. We treat separately the first and second types~\eqref{eq:f_type1} and~\eqref{eq:f_type2}.

\begin{theorem}[Malthusian growth: blow up vs extinction]\label{th:extinction_vs_explo:lambda}
Assume that $f_1,f_2$ are of the first type~\eqref{eq:f_type1}, let $\lambda$ be given by \eqref{eq:def_lambda} for the operator $\A$ defined in~\eqref{eq:defAbis} with fitnesses~\eqref{def:ribis}. Let $\densbf$ be the solution of~\eqref{eq:sys_main1bis} given by Theorem~$\ref{thm:existunidens}$-{\rm{(i)}}, with initial condition $\densbf^0=(u^0_1,u^0_2)$ satisfying~{\rm{(H1)-(H3)}}, and let $N_1(t)$ and $N_2(t)$ be its population sizes in each habitat.
\begin{itemize}
\item[{\rm{(i)}}] If $\lambda<0$, then $N_1(t) \to +\infty$ and $N_2(t)\to+\infty$ as $t\to +\infty$ $($blow up of the population$)$.
\item[{\rm{(ii)}}] If $\lambda=0$ and if $u^0_1$ and $u^0_2$ are compactly supported, then:
\[\limsup_{t\to+\infty}\big(N_1(t)+N_2(t)\big)<+\infty\quad \text{$($boundedness of the population$)$.}\]
Furthermore, there exist bounded positive stationary solutions of~\eqref{eq:sys_main1bis} with finite population sizes.
\item[{\rm{(iii)}}] If $\lambda>0$ and if $u^0_1$ and $u^0_2$ are compactly supported, then $N_1(t)+N_2(t)\to 0$ as $t\to +\infty$ $($extinction of the population$)$.
\end{itemize}
\end{theorem}

\begin{theorem}[Logistic growth: persistence vs extinction]\label{th:extinction_vs_explo:logistic}
Assume that $f_1,f_2$ are of the second type~\eqref{eq:f_type2}, and let $\lambda$ given by~\eqref{eq:def_lambda} for the operator $\A$ defined in~\eqref{eq:defA} with fitnesses~\eqref{def:ri}. Let $\densbf$ be the solution of~\eqref{eq:sys_main1} given by Theorem~$\ref{thm:existunidens}$-{\rm{(ii)}}, with initial condition $\densbf^0=(u^0,u^0\circ\iota)$ satisfying~\eqref{hyp:HS} and~{\rm{(H1)-(H3)}}, and let $N(t) = N_1(t)=N_2(t)$ be its population size in each habitat.
\begin{itemize}
\item[{\rm{(i)}}] If $\lambda<0$, then:
\be\label{persistence0}
\displaystyle \limsup_{t\to+\infty}N(t)>0 \quad \text{$($persistence of the population$)$},
\ee
and even:
\be\label{persistence}
\displaystyle 0\!<\!\liminf_{t\to+\infty}N(t)\!\le\!\limsup_{t\to+\infty}N(t)\!<\!+\infty\ \text{$($persistence of the population, in a strong sense$)$}
\ee
for some initial conditions~$\densbf^0$ $($see Remark~$\ref{rempersistence}$ below$)$.
\item[{\rm{(ii)}}] If $\lambda\ge 0$ and if  $u^0$ is compactly supported, then $N(t) \to 0$ as $t\to +\infty$ $($extinction of the population$)$.
\end{itemize}
\end{theorem}

As a consequence of Theorems~\ref{th:extinction_vs_explo:lambda}-\ref{th:extinction_vs_explo:logistic}, the fate of the population is determined by the sign of~$\lambda$, i.e., by the linear stability of the steady state $\densbf=(0,0)$, whether the growth functions~$f_i$ be of the first or second type. The main differences between the Malthusian case and the logistic case arise when this steady state $(0,0)$ is unstable ($\lambda<0$). Although persistence occurs with both types of growth functions, the population size remains bounded with type 2 growth functions, due to the nonlocal competition term. We conjecture that it converges to $-\lambda$, as $t\to +\infty$. Interestingly, the threshold case $\lambda=0$ leads to very different behaviours, depending on the type of growth functions: in the absence of competition (Malthusian growth), persistence is still possible in this case, although it is not in the logistic case. Biologically, however, the particular case $\lambda=0$ is presumably not relevant. 

\begin{remark}\label{rempersistence}{\rm In part~(i) of Theorem~$\ref{th:extinction_vs_explo:logistic}$, namely if $\lambda\!<\!0$, the initial conditions $\densbf^0\!=\!(u^0,u^0\circ\iota)$ such that~\eqref{persistence} holds are those which are trapped between two multiples of the principal eigenfunctions associated to the operator $\A$ given by~\eqref{eq:defA} but acting this time on $[W^{2,n}_{loc}(\R^n)\cap C_0(\R^n)]^2$, where~$C_0(\R^n)$ is the space of continuous functions in $\R^n$ converging to $0$ at infinity. Such eigenfunctions are introduced in Lemma~\ref{lem:fpp} below.}
\end{remark}

In the following results, for problem~\eqref{eq:sys_main1} and the operator $\A$ defined in~\eqref{eq:defA} with fitnesses~\eqref{def:ri}, we now use the above persistence/extinction criteria to study the effect of the parameters, namely the migration rate~$\delta>0$, the habitat difference $m_D\ge0$ given by~\eqref{def:mD}, the mutational parameter $\mu>0$ and the fitness optimum~$\rmax$, on the fate of the population. 

\begin{prop}\label{prop:continuity_lambda}
For the operator $\A$ defined in~\eqref{eq:defA} with fitnesses~\eqref{def:ri}, let $\lambda_{\delta,m_D,\mu,\rmax}:=\lambda$ be given by~\eqref{eq:def_lambda}. Then the map:
$$(\delta,m_D,\mu,\rmax)\mapsto\lambda_{\delta,m_D,\mu,\rmax},$$
is continuous in $(0,+\infty)\times[0,+\infty)\times(0,+\infty)\times\R$. Moreover,
\begin{itemize}
\item[(i)] for each $(m_D,\mu,\rmax)\in[0,+\infty)\times(0,+\infty)\times\R$, the map $\delta\mapsto\lambda_{\delta,m_D,\mu,\rmax}$ is nondecreasing $($and even increasing if $m_D>0)$ and concave in $(0,+\infty)$, and:
$$\lim_{\delta \to 0^+} \lambda_{\delta,m_D,\mu,\rmax}=-\rmax +{\dd\,n \over 2},\quad\lim_{\delta \to +\infty} \lambda_{\delta,m_D,\mu,\rmax} =-\rmax +{\dd\,n \over 2}+{m_D\over 4};$$
\item[(ii)] for each $(\delta,\mu,\rmax)\in(0,+\infty)\times(0,+\infty)\times\R$, the map $m_D\mapsto\lambda_{\delta,m_D,\mu,\rmax}$ is increasing in~$[0,+\infty)$
and:
$$\lambda_{\delta,0,\mu,\rmax}=-\rmax +{\dd\,n \over 2},\quad\lim_{m_D\to +\infty} \lambda_{\delta,m_D,\mu,\rmax} =-\rmax +{\dd\,n \over 2}+\delta;$$
\item[(iii)] for each $(\delta,m_D,\rmax)\in(0,+\infty)\times[0,+\infty)\times\R$, the map $\mu\mapsto\lambda_{\delta,m_D,\mu,\rmax}$ is increasing in~$(0,+\infty)$, and there is $\lambda^0\in[0,\min(\delta,m_D/4)]$, independent of $\rmax$, such that:
$$\lim_{\mu\to 0^+} \lambda_{\delta,m_D,\mu,\rmax}=-\rmax+\lambda^0,\quad\lim_{\mu\to +\infty} \lambda_{\delta,m_D,\mu,\rmax} =+\infty;$$
\item[(iv)] for each $(\delta,m_D,\mu)\in(0,+\infty)\times[0,+\infty)\times(0,+\infty)$, one has $\lambda_{\delta,m_D,\mu,\rmax}=-\rmax+\lambda_{\delta,m_D,\mu,0}$ for all $\rmax\in\R$, hence $\lambda_{\delta,m_D,\mu,\rmax}\to\mp\infty$ as $\rmax\to\pm\infty$.
\end{itemize}
\end{prop}

A corollary of Theorems~\ref{th:extinction_vs_explo:lambda}-\ref{th:extinction_vs_explo:logistic} and Proposition~\ref{prop:continuity_lambda} follows immediately with straightforward proof.

\begin{cor}\label{cor:explo_extinc_dens}
For problem~\eqref{eq:sys_main1} with fitnesses~\eqref{def:ri}, let $\densbf$ be the solution given by Theorem~$\ref{thm:existunidens}$-{\rm{(ii)}}, with initial condition $\densbf^0=(u^0,u^0\circ\iota)$ satisfying~\eqref{hyp:HS} and~{\rm{(H1)-(H3)}}, and let $N(t) = N_1(t)=N_2(t)$ be its population size in each habitat.
\begin{itemize}
\item[{\rm(i)}] If  $\rmax\ge\dd\,n/2+\min(m_D/4,\delta)$ and $m_D>0$, then $\lim_{t\to+\infty}N(t)=+\infty$ for the first type~\eqref{eq:f_type1}, whereas~\eqref{persistence0} holds and even~\eqref{persistence} is satisfied for some initial conditions $\densbf^0$ for the second type~\eqref{eq:f_type2}.
\item[{\rm(ii)}] If $\dd\,n/2<\rmax<\dd\,n/2+m_D/4$, then there exists $\delta_{crit}>\rmax-\dd\,n/2$, independent of $\densbf^0$, such that:
\begin{itemize}
\item if $\delta<\delta_{crit}$, then $\lim_{t\to+\infty}N(t)=+\infty$ for the first type~\eqref{eq:f_type1}, whereas~\eqref{persistence0} holds and even~\eqref{persistence} is satisfied for some initial conditions $\densbf^0$ for the second type~\eqref{eq:f_type2};
\item if $\delta=\delta_{crit}$ and if $u^0$ is compactly supported, then $\limsup_{t\to+\infty}N(t)<+\infty$ for the first type~\eqref{eq:f_type1} and $N(t)\to0$ as $t\to+\infty$ for the second type~\eqref{eq:f_type2};
\item if $\delta>\delta_{crit}$ and if $u^0$ is compactly supported, then $\lim_{t\to+\infty}N(t)=0$ for both types~\eqref{eq:f_type1} and~\eqref{eq:f_type2}.
\end{itemize}
\item[{\rm(iii)}] If $\dd\,n/2<\rmax<\dd\,n/2+\delta$, then there exists $m_{D,crit}>4\,(\rmax-\dd\,n/2)$, independent of $\densbf^0$, such that:
\begin{itemize}
\item if $m_D<m_{D,crit}$, then $\lim_{t\to+\infty}N(t)=+\infty$ for the first type~\eqref{eq:f_type1}, whereas~\eqref{persistence0} holds and even~\eqref{persistence} is satisfied for some initial conditions $\densbf^0$ for the second type~\eqref{eq:f_type2};
\item if $m_D=m_{D,crit}$ and if $u^0$ is compactly supported, then $\limsup_{t\to+\infty}N(t)\!<+\infty$ for the first type~\eqref{eq:f_type1} and $N(t)\to0$ as $t\to+\infty$ for the second type~\eqref{eq:f_type2};
\item if $m_D>m_{D,crit}$ and if $u^0$ is compactly supported, then $\lim_{t\to+\infty}N(t)=0$ for both types~\eqref{eq:f_type1} and~\eqref{eq:f_type2}.
\end{itemize}
\item[{\rm(iv)}] If $\rmax\le\max({\dd\,n}/2,\lambda^0)$ (with $\lambda^0$ defined in Proposition~\ref{prop:continuity_lambda}) and if $u^0$ is compactly supported, then $\lim_{t\to+\infty}N(t)=0$ for both types~\eqref{eq:f_type1} and~\eqref{eq:f_type2}.
\item[{\rm(v)}] If $\rmax>\lambda^0$ (with $\lambda^0$ defined in Proposition~\ref{prop:continuity_lambda}), then  there exists $\dd_{crit}>0$, independent of $\densbf^0$, such that $\dd_{crit}>(2/n)\times(\rmax-\min(\delta,m_D/4))$ if $m_D>0$ $($resp. $\dd_{crit}=(2/n)\times\rmax$ if $m_D=0$$)$ and:
\begin{itemize}
\item if $\dd<\dd_{crit}$, then $\lim_{t\to+\infty}N(t)=+\infty$ for the first type~\eqref{eq:f_type1}, whereas~\eqref{persistence0} holds and even~\eqref{persistence} is satisfied for some initial conditions $\densbf^0$ for the second type~\eqref{eq:f_type2};
\item if $\dd=\dd_{crit}$ and if $u^0$ is compactly supported, then $\limsup_{t\to+\infty}N(t)<+\infty$ for the first type~\eqref{eq:f_type1} and $N(t)\to0$ as $t\to+\infty$ for the second type~\eqref{eq:f_type2};
\item if $\dd>\dd_{crit}$ and if $u^0$ is compactly supported, then $\lim_{t\to+\infty}N(t)=0$ for both types~\eqref{eq:f_type1} and~\eqref{eq:f_type2}.
\end{itemize}
\item[{\rm(vi)}] For every $\delta>0$, $m_D\ge0$ and $\dd>0$, there exists $r_{\max,crit}\in\R$, independent of $\densbf^0$, such that $\max(\dd\,n/2,\lambda^0)<r_{\max,crit}<\dd\,n/2+\min(\delta,m_D/4)$ if $m_D>0$ $($resp. $r_{\max,crit}=\dd\,n/2$ if $m_D=0$$)$, and:
\begin{itemize}
\item if $\rmax<r_{\max,crit}$ and if $u^0$ is compactly supported, then $\lim_{t\to+\infty}N(t)=0$ for both types~\eqref{eq:f_type1} and~\eqref{eq:f_type2};
\item if $\rmax=r_{\max,crit}$ and if $u^0$ is compactly supported, then $\limsup_{t\to+\infty}N(t)<+\infty$ for the first type~\eqref{eq:f_type1} and $N(t)\to0$ as $t\to+\infty$ for the second type~\eqref{eq:f_type2};
\item if $\rmax>r_{\max,crit}$, then $\lim_{t\to+\infty}N(t)=+\infty$ for the first type~\eqref{eq:f_type1}, whereas~\eqref{persistence0} holds and even~\eqref{persistence} is satisfied for some initial conditions $\densbf^0$ for the second type~\eqref{eq:f_type2}.
\end{itemize}
\end{itemize}
\end{cor}

An interpretation of parts~(i)-(iii) of Corollary~\ref{cor:explo_extinc_dens} is that, when the fitness optimum $\rmax$ is larger than the mutation load, namely~$\rmax>\dd\,n/2$, and the habitat difference $m_D$ is low, namely $m_D\leq4\,(\rmax-\dd\,n/2)$, the population can adapt to the global environment, whatever the migration rate~$\delta$. However, when the habitat difference is high, namely $m_D>4(\rmax-\dd\,n/2)$, the population can only survive if the migration rate is low~($\delta \le \delta_{crit}$). Conversely, under the same condition on $\rmax$, when the migration rate $\delta$ is low, namely $\delta\leq\rmax-\dd\,n/2$, the population can adapt to the global environment, whatever the habitat difference $m_D$. However, when the migration rate is high, namely $\delta>\rmax-\dd\,n/2$, the population can only survive if the habitat difference is low ($m_D\le m_{D,crit}$). These results are coherent with the numerical simulations of Figure~\ref{fig:graph_mDcrit}. See Section~\ref{sec:discu} for more detailed interpretations of these results. 

The last result is related to the generalist nature of the population in the limit of infinite migration rates $\delta$.

\begin{theorem}\label{thm:same_dens_large_delta}
For problem~\eqref{eq:sys_main1} with fitnesses~\eqref{def:ri} and growth functions $f_1,f_2$ of the first type~\eqref{eq:f_type1}, let $\densbf_\delta=(\dens_{\delta,1},\dens_{\delta,2})$ be the solution given by Theorem~$\ref{thm:existunidens}$-{\rm{(ii)}}, with a fixed initial condition $\densbf^0=(u^0_1,u^0_2)$, independent of $\delta$, satisfying the assumptions~{\rm{(H1)-(H3)}} $($but the assumption~\eqref{hyp:HS} may not be satisfied$)$. Then:
\[
\lim_{\delta \to +\infty}\|\dens_{\delta,1}(t,\cdot)-\dens_{\delta,2}(t,\cdot)\|_{L^\infty(\R^n)}=0,\ \hbox{ locally uniformly in }t\in(0,+\infty).
\]
\end{theorem}

In other words, a strong migration rate $\delta$ merges the two populations into one global population, since the exchanges between them are very large. The population then goes to be generalist at every time $t>0$, even if it is not initially. This is consistent with current knowledge, see \cite{DebRon13}.

\section{Discussion \label{sec:discu}}

\paragraph{On the biological interpretation of the main results.} 
For problem~\eqref{eq:sys_main1} with fitnesses~\eqref{def:ri}, Proposition~\ref{prop:continuity_lambda} together with Theorems~\ref{th:extinction_vs_explo:lambda}-\ref{th:extinction_vs_explo:logistic} show that the more the two environments are connected by migration (i.e., when $\delta$ is increased), the lower are the chances of persistence. In the absence of migration, when the two habitats are not connected ($\delta=0$), it was already known that persistence occurs if  $\rmax > \dd\,n /2$ \cite{HamLavMarRoq20,MarRoq16}, whereas $\rmax < \dd\,n /2$ leads to extinction (for both types of growth functions). In the case $\delta=0$, at large times, the mean fitness $\rb(t)$ converges to $\rmax - \dd\,n /2$,  with $\dd\, n/2$ the mutation load. As already mentioned, if the mutation load exceeds $\rmax$, the population is doomed to extinction. 

When $\delta$ becomes positive, some individuals migrate between the two environments. Gene\-rally these individuals are better adapted to their environment of origin. Thus, as shown by Proposition~\ref{prop:continuity_lambda}, increasing the migration rate increases the global maladaptation. Ultimately, when $\delta\to +\infty$, the condition for persistence becomes $\rmax>m_D/4+\dd\,n/2$: in this case, as shown by Theorem~\ref{thm:same_dens_large_delta}, the two phenotypic populations merge into a single one, centered at the origin, in-between the two optima. We observe that in addition to the mutation load $\dd \, n/2$, a \emph{migration load} $m_D/4$ appears. It is proportional to the habitat difference $m_D$.

Increasing the habitat difference $m_D$ also increases global maladaptation (Proposition~\ref{prop:continuity_lambda}, part~ii). Although this result seems natural, its mathematical proof is rather involved. It uses the fact that the eigenfunctions associated with each optimum are asymmetric and biased towards the side of the other optimum, see Section~\ref{sec:33} and in particular Eq.~\eqref{claimpsi}. This asymmetry reflects the advantage, for a phenotype $\x$ in a given habitat $i$ and with a fixed fitness $r_i(\x)$, to be closer to the phenotype optimum of the other habitat while keeping the same fitness value. When the habitat difference becomes very large ($m_D \to +\infty$), the condition for persistence becomes $\rmax>\delta+\dd\,n/2$ which means that the migration rate $\delta$ plays the same role as a death rate, and the migration load is simply $\delta$.

If $\dd\,n/2<\rmax<m_D/4+\dd\,n/2$, populations are doomed to extinction for large migration rates, but survive for small migration rates. Corollary~\ref{cor:explo_extinc_dens} shows that there exists a migration threshold such that persistence is possible if the migration rate is below this threshold, but not if the migration rate is above this threshold. Thus, increasing the migration rate may imply a `lethal migration effect', comparable to lethal mutagenesis. If $\rmax\ge m_D/4+\dd\,n/2$, persistence always occurs, independently of the migration rate.

Similarly, Corollary~\ref{cor:explo_extinc_dens} shows that if $\dd\,n/2<\rmax<\delta+\dd\,n/2$, populations go extinct for large habitat differences and persist for small habitat differences, with again a threshold value of $m_D$ which determines persistence. If $\rmax\ge \delta+\dd\,n/2$, increasing the habitat difference will have no effect on the persistence of the population.  

Mutation also has a detrimental effect on persistence. As in the case of a single habitat, increasing the mutation term $\mu$ ultimately leads to lethal mutagenesis (Proposition~\ref{prop:continuity_lambda}, part iii). However, contrarily to the effects observed above for $\delta$ and $m_D$, this phenomenon occurs whatever the value of $\rmax.$

\paragraph{Implications in agroecology.} One of the fundamental principles in agroecology is to promote diversified agroecosystems rather than uniform cultures \cite{CaqGas20,FAO18}. Some empirical study already illustrated the higher resilience of such diversified agroecosystems \cite{BorLar18} to plant diseases. In our case, the two environments can be interpreted as two different types of host plants (different species, or different genetic variants) and the populations of phenotypes $u_1$, $u_2$ describe the density of a pathogen over these two types of host plants. With this interpretation, our study advocates for more diversified cultures, with strong migration of the pathogens between the host plants: it should reduce the chances of persistence of the pathogen over the agroecosystem. This is consistent, therefore, with the above-mentioned principle of plant diversification. However, we point out that this conclusion may not be valid for three environments or more: as discussed in \cite{LavMar19}, the presence of a third environment associated with a phenotype optimum  between the two others may lead to higher chances of persistence of the pathogen, compared to two environments, due to a `springboard' effect. By now, and up to our knowledge, there is no rigorous mathematical proof of this result. 

\paragraph{On the derivation of quantitative estimates.} The methods used in our paper do not allow for a computation of the migration load: when $\delta=0$ in~\eqref{eq:sys_main1} with fitnesses~\eqref{def:ri}, as discussed above, the mean growth rate $\rb(t)$ converges to $\rmax-\mu\, n/2$. With positive values of $\delta$, it should converge to some value  $\rmax-\mu\, n/2- \text{Load}_{migr}(\delta)$, with $\text{Load}_{migr}(\delta) \in (0,\min(\delta,m_D/4))$, the migration load. The determination of $\text{Load}_{migr}(\delta)$ would help disentangling the respective effect of mutation and migration on the persistence of a population. Additionally, Theorem~\ref{thm:same_dens_large_delta} shows that when the migration rate is increased the two population merge into a single one, which may be qualified as `generalist'. This is consistent with the results that have been obtained by \cite{MirGan20} in the case $n=1$ with methods based on constrained Hamilton-Jacobi equations, and more broadly with current literature \cite{DebRon13,MesCzi97,PapDav13} in evolutionary biology. This means that the mean phenotype in each environment converges to $\x=0$. With smaller migration rates, the two populations should behave as `specialists',  with mean phenotypes that converges to $\opt_1$ and $\opt_2$ respectively as $\delta\to 0$. In a forthcoming work, using the methods in \cite{HamLavMarRoq20} based on the analysis of moment generating functions associated with the distribution of fitness, we will aim to derive quantitative estimates for the migration load, the lethal migration threshold $\delta_{crit}$ and the respective distributions of phenotypes in the two environments.

\section{Proofs}\label{section:proofs}

This section is devoted to the proofs of the results stated in Section~\ref{sec:results}. Section~\ref{sec:31} is devoted to the proof of Theorem~\ref{thm:existunidens} on the well-posedness of the Cauchy problems~\eqref{eq:sys_main1} and~\eqref{eq:sys_main1bis}. Section~\ref{sec:32} is concerned with the proof of Theorems~\ref{th:extinction_vs_explo:lambda}-\ref{th:extinction_vs_explo:logistic} on the large time behaviour of the population size, and Section~\ref{sec:33} with the dependence of the fate of the population with respect to the parameters.

\subsection{The Cauchy problems~\eqref{eq:sys_main1} and~\eqref{eq:sys_main1bis}}\label{sec:31}

\noindent \textit{Proof of Theorem~\ref{thm:existunidens}}.
Part (i). We begin with problem~\eqref{eq:sys_main1bis} associated to fitnesses~\eqref{def:ribis}, growth functions $f_1,f_2$ of the first type~\eqref{eq:f_type1}, and initial conditions $\densbf^0=(\dens_1^0,\dens^0_2)$ satisfying~(H1)-(H3). As we will see later in part~(ii) of the proof, the results in the case of problem~\eqref{eq:sys_main1} with fitnesses~\eqref{def:ri} and growth functions $f_1,f_2$ of the second type~\eqref{eq:f_type2} are then straightforward thanks to a change of functions when the initial conditions $\densbf^0=(\dens_1^0,\dens^0_2)$ satisfy~(H1)-(H3) together with~\eqref{hyp:HS}.

So, let us first assume that $f_1,f_2$ are of the first type~\eqref{eq:f_type1}. Thanks to the assumptions~(H1)-(H3), owing to the definition~\eqref{def:ribis} of the fitnesses $r_i$ and setting $R_{\max}=\max(r_{\max,1},r_{\max,2})$, it follows from~\cite[Theorem~3]{Bes79} that, for any $T>0$, the Cauchy problem:
\[\left\{\begin{array}{rcll}
\partial_tv_1 (t,\x) & = &\ds \frac{\dd^2} 2 \ \Delta v_1(t,\x)+[r_1(\x)-R_{\max}]\,v_1(t,\x)-\delta_{1,1}v_1(t,\x)+\delta_{1,2}v_2(t,\x),   & t\ge0, \ \x\in \R^n, \vspace{1mm}\\
     \partial_t v_2 (t,\x) & =& \ds\frac{\dd^2} 2 \ \Delta v_2(t,\x)+[r_2(\x)-R_{\max}]\, v_2(t,\x) +\delta_{2,1}v_1(t,\x)-\delta_{2,2}v_2(t,\x),   & t\ge0, \ \x\in \R^n, \vspace{1mm}\\
     \mathbf v(0, \x) & = & \densbf^0(\x)=(u^0_1(\x),u^0_2(\x)),  \  \x \in \R^n, &
\end{array}\right.\]
admits a solution $\mathbf v=(v_1,v_2)\in[C^{1,2}([0,T] \times \R^n)\cap L^\infty((0,T)\times\R^n)]^2$, such that $\mathbf v(t,\x)\to(0,0)$ as~$\|\x\|\to+\infty$ uniformly in $t\in[0,T]$. Thus, the function $\densbf:(t,\x)\mapsto e^{R_{\max}t}\, \mathbf v(t,\x)$, defined in~$[0,T]\times \R^n$, is a bounded classical solution of~\eqref{eq:sys_main1bis} satisfying the same properties as~$\mathbf v$. Moreover, this solution is nonnegative (componentwise) from the comparison principle~\cite[Lemma~1]{Wei75} applied to this linear cooperative system. This maximum principle also yields the uniqueness of this solution $\densbf$. Since the initial population density in each habitat is not identically equal to~$0$ by assumption~(H2), the nonnegativity of each component $\dens_i$ and the strong parabolic maximum principle applied to each linear operator $\partial_t-(\dd^2/ 2)\Delta -r_i(\x)+\delta_{i,i}$ (for $i\in\{1,2\}$) yield the posi\-tivity of each component $\dens_i$ in $(0,T]\times\R^n$. As $T>0$ can be chosen arbitrarily, these existence, uniqueness and positivity results extend to $t\in (0,+\infty)$, with local boundedness in $t$.

For part~(i) of the proof, it still remains to show that the population sizes and mean fitnesses $N_i(t)$ and $\rb_i(t)$ defined by~\eqref{dfn:pop_size_pde}-\eqref{dfn:meanfit} are real valued, continuous and satisfy~\eqref{eq:ODEmalthusbis}. We first establish some bounds and, to do so, we construct a super-solution for  $\densbf=(\dens_1,\dens_2)$. Let us first denote:
$$\omega\!=\!r_{\max,2}\!-\!r_{\max,1}\!+\!\delta_{1,1}\!-\!\delta_{2,2},\ \ \rho\!=\!\frac{r_{\max,1}\!+\!r_{\max,2}}{2}\!-\!\frac{\delta_{1,1}\!+\!\delta_{2,2}}{2}\ \hbox{ and }\ \gamma\!=\!\frac{\sqrt{\omega^2\!+\!4\delta_{1,2}\delta_{2,1}}}{2}>0,$$
and let us set, for all $t>0$ and~$\x\in\R^n$:
\begin{multline}
\mathbf h (t,\x):= \left( \begin{matrix} h_1(t,\x) \\ h_2(t,\x)\end{matrix}\right)  :=  \displaystyle e^{\rho t} \left[K_{t} \ast \dens^0_1 \right](\x) \left( \begin{matrix} \displaystyle \cosh(\gamma t)-\frac{\omega\sinh(\gamma t)}{2\gamma} \\ \displaystyle\frac{\delta_{2,1}\sinh(\gamma t)}{\gamma}\end{matrix} \right)\\
 \displaystyle+\,e^{\rho t} \left[K_{t} \ast \dens^0_2 \right](\x) \left( \begin{matrix} \displaystyle\frac{\delta_{1,2}\sinh(\gamma t)}{\gamma} \\ \displaystyle\cosh(\gamma t)+\frac{\omega\sinh(\gamma t)}{2\gamma}\end{matrix} \right),
\label{dfn:hbf}
\end{multline}
with:
\[
\forall\,t>0, \ \forall\,\x \in \R^n, \ \ K_{t}(\x) = \frac{e^{-\|\x\|^2/(2\dd^2 t)}}{(2\pi  \dd ^2 t) ^{n/2}},
\]
and $\mathbf h(0,\x) = \densbf^0(\x)=\densbf(0,\x)$ (we use in the explicit expression~\eqref{dfn:hbf} the fact that the mutational parameter $\mu$ is the same in both habitats). The function $\mathbf h$ is of class $[C^\infty((0,+\infty)\times\R^n)\cap C([0,+\infty)\times\R^n)]^2$, it is locally bounded in time, it converges to $(0,0)$ as $\|\x\|\to+\infty$ locally uniformly in $t\in\R_+$, and it is straightforward to check that it satisfies:
\[
\partial_ t \mathbf h(t,\x) = {\mu^2\over 2}\Delta \mathbf h(t,\x) +\left(\begin{matrix}
r_{\max,1}-\delta_{1,1} & \delta_{1,2} \\
\delta_{2,1} & r_{\max,2} -\delta_{2,2}
\end{matrix}\right) \mathbf h (t,\x),
\]
for all $t>0$ and $\x \in \R^n$. Let $\psibf(t,\x) := \densbf(t,\x) -\mathbf h(t,\x)$. We see that $\psibf(0,\x)=(0,0)$ for all~$\x\in\R^n$, and:
\begin{equation}\label{eqpsi}
\partial_t \psibf(t,\x) - {\mu^2\over 2}\Delta \psibf(t,\x) -\left(\begin{matrix}
r_{\max,1}-\delta_{1,1} & \delta_{1,2} \\
\delta_{2,1} & r_{\max,2} -\delta_{2,2}
\end{matrix}\right) \psibf(t,\x) = \left( \begin{matrix} m_1(\x)\,\dens_1(t,\x) \\
m_2(\x)\,\dens_2(t,\x)
\end{matrix} \right)\le\left(\begin{matrix} 0\\ 0\end{matrix}\right),
\end{equation}
for all $t>0$ and $\x \in \R^n$, with:
\be\label{def:m(x)}
m_i(\x):=r_i(\x)-r_{\max,i}=-\frac{\|\x-\opt_i\|^2}{2}\le0.
\ee
Again, the comparison principle~\cite[Lemma~1]{Wei75}  implies that $\psibf\le 0$ (componentwise) in $\R_+\times \R^n$, hence:
\be\label{eq:maj1:dens}
0\le\densbf(t,\x) \le \mathbf h(t,\x)\ \hbox{ for all }(t,\x)\in\R_+\times\R^n.
\ee 
The strong parabolic maximum principle actually implies that the second inequality, as is the first one, is strict in $(0,+\infty)\times\R^n$, as follows from~\eqref{eqpsi} together with the positivity of $\dens_1$ and~$\dens_2$. Moreover, for $i\in\{1,2\}$ and $t>0$,
$$\int_{\R^n} \left[K_{t} \ast \dens^0_i \right](\x)\,\md\x =\int_{\R^n}u_i^0(\x)\,\md\x=:N_i^0 <+\infty.$$
Thus, $\int_{\R^n}h_i(t,\x)\,\md\x\le (N_1^0+N_2^0)\,M(t)$ for all $t>0$, 
with:
\[M(t):=e^{\rho t}\times\max\Big(\cosh(\gamma t)+\frac{|\omega|\sinh(\gamma t)}{2\gamma},\frac{\max(\delta_{1,2},\delta_{2,1})\,\sinh(\gamma t)}{\gamma}\Big).\]
Hence, from~\eqref{eq:maj1:dens} and the positivity of~$\dens_i$ in~$(0,+\infty)\times\R^n$, there holds:
\be\label{Nit}
0<N_i(t)=\int_{\R^n}u_i(t,\x)\,\md\x\le(N_1^0+N_2^0)\,M(t),
\ee
for all $t>0$, as well as for $t=0$ trivially.

Consider now any time $t\ge0$ and let us prove that $\rb_i(t)$ defined in~\eqref{dfn:meanfit} is finite, for $i\in\{1,2\}$. First, the hypotheses (H2)-(H3) imply that $\rb_i(0)$ is finite. Assume then that $t>0$. From~\eqref{def:m(x)}-\eqref{eq:maj1:dens} and the positivity of $\dens_i$, we have:
\be\label{ineq:Nit}
r_{\max,i}\, N_i(t) \ge \int_{\R^n} r_i(\x) \, \dens_i(t,\x)\,\md\x \ge r_{\max,i}\, N_i(t)+ \int_{\R^n} m_i(\x) \, h_i(t,\x)\,\md\x.
\ee
Thus, to show that $\rb_i(t)$ is finite, we only have to show that the last term in the right-hand side of the above equation is finite. First, we note that:
\begin{equation} \label{eq:ineghi}
0\le h_i(t,\cdot) \le M(t)\,K_{t} \ast (\dens^0_1+\dens^0_2)\ \hbox{ in }\R^n.
\end{equation}
Then, still using the assumption (H3), we have:
\begin{align*}
0 & \le\!\!\int_{\R^n}\!\!-m_i(\x)\,[K_{t} \ast  \dens_i^0] (\x) \, \md\x , \\
& \le  {1 \over 2(2\pi \dd ^2 t) ^{n/2}}\int_{\R^n} \int_{\R^n}  \ \|\x-\opt_i\|^2 \ e^{-\|\x-\y\|^2/(2\dd^2 t)}  g(\|\y-\opt_i\|)\, \md\y \, \md\x, \\
& ={1\over 2\pi  ^{n/2}}\int_{\R^n} \int_{\R^n}  \ \|\x-\opt_i\|^2 \ e^{-\|\z\|^2}  g(\|\x-\dd \sqrt{2t} \,\z-\opt_i\|)\, \md\z \, \md\x, \\
& \le  {1\over 2\pi  ^{n/2}}\int_{\R^n} \int_{\|\z\|\le \|\x-\opt_i\|/(2\dd\sqrt{2t})}  \ \|\x-\opt_i\|^2 \ e^{-\|\z\|^2}  g\left({\|\x-\opt_i\|\over 2}\right)\, \md\z \, \md\x \\
& \qquad \qquad \qquad \qquad\quad \qquad   +{1 \over 2\pi  ^{n/2}}\int_{\R^n} \int_{\|\z\|>\|\x-\opt_i\|/(2\dd\sqrt{2t})}  \ \|\x-\opt_i\|^2 \ e^{-\|\z\|^2}  g(0)\, \md\z \, \md\x, \\
& \le {1\over 2\pi ^{n/2}} \left[ \pi ^{n/2} \int_{\R^n} \|\x-\opt_i\|^2 \ g\left({\|\x-\opt_i\|\over 2}\right)\, \md\x+ g(0) \int_{\R^n} \zeta_{t}(\|\x-\opt_i\|) \|\x-\opt_i\|^2 \, \md\x\right],
\end{align*}
\normalsize
where:
\be\label{defzeta}
\zeta_{t}(r):=\int_{\|\z\|\ge r/(2{\dd}\sqrt{2t})}e^{-\|\z\|^2}\,\md\z=O\big(e^{-r}\big),\ \hbox{ as }r\to+\infty.
\ee
The assumption (H3) thus implies that:
\begin{equation} \label{eq:K*1}
    0\le \int_{\R^n} -m_i(\x) [K_{t} \ast  \dens_i^0] (\x)\,\md\x<+\infty,
\end{equation}
for every $t>0$. Let us now check that $-\int_{\R^n} m_i(\x) [K_{t} \ast  \dens_j^0] (\x) \, \md\x<+\infty$ for $i\neq j\in \{1,2\}$:
\begin{align*}
\int_{\R^n} -m_i(\x) & [K_{t} \ast  \dens_j^0] (\x)\, \md\x  \\ & \le  {1 \over 2(2\pi \dd ^2 t) ^{n/2}}\int_{\R^n} \int_{\R^n}  \ \|\x-\opt_i\|^2 \ e^{-\|\x-\y\|^2/(2\dd^2 t)}  g(\|\y-\opt_j\|)\, \md\y \, \md\x, \\
& ={1 \over 2\pi  ^{n/2}}\int_{\R^n} \int_{\R^n}  \, \|\x-\opt_i\|^2 \ e^{-\|\z\|^2}  g(\|\x-\dd \sqrt{2t} \,\z-\opt_j\|)\, \md\z \, \md\x, \\
& \le  {1 \over 2\pi  ^{n/2}}\int_{\R^n} \int_{\|\z\|\le \|\x-\opt_j\|/(2\dd\sqrt{2t})}  \ \|\x-\opt_i\|^2 \ e^{-\|\z\|^2}  g\left({\|\x-\opt_j\|\over 2}\right)\, \md\z \, \md\x \\
& \quad\,\qquad \qquad ^  +{1 \over 2\pi  ^{n/2}}\int_{\R^n} \int_{\|\z\|>\|\x-\opt_j\|/(2\dd\sqrt{2t})}  \ \|\x-\opt_i\|^2 \ e^{-\|\z\|^2}  g(0)\, \md\z \, \md\x, \\
& \le {1 \over 2\pi  ^{n/2}} \left[ \pi ^{n/2} \int_{\R^n}\big(2\|\x-\opt_j\|^2+8\beta^2\big)\,g\left({\|\x-\opt_j\|\over 2}\right)\,  \md\x  \right.\\
& \qquad \qquad \qquad \qquad \quad\quad\, \left. +g(0) \int_{\R^n} \zeta_{t}(\|\x-\opt_j\|)\,\big(2\|\x-\opt_j\|^2+8\beta^2 \big)\, \md\x\right],
\end{align*}\normalsize
where we recall that $\beta$ is defined in~\eqref{def:opt_beta}. Thus,  (H3) implies that:
\begin{equation} \label{eq:K*2}
    0\le \int_{\R^n} -m_i(\x) [K_{t} \ast  \dens_j^0] (\x)\,  \md\x<+\infty.
\end{equation}
Adding \eqref{eq:K*1} and \eqref{eq:K*2}, and using~\eqref{eq:ineghi}, we obtain that:
$$0\le\int_{\R^n} -m_i(\x) \,h_i (t,\x)\,  \md\x<+\infty,$$
and, together with~\eqref{ineq:Nit}, we infer that $-\infty<\rb_i(t)\le r_{\max,i}$ for $i\in\{1,2\}$ and  $t>0$ (and also for $t=0$ as already emphasized). 

Finally, since the quantities $\zeta_{t}(r)$ given in~\eqref{defzeta} are nondecreasing with respect to $t>0$, the same arguments as above together with Lebesgue's dominated convergence theorem yield the continuity of the maps $t\mapsto N_i(t)$, $t\mapsto\int_{\R^n}m_i(\x)\,u_i(t,\x)\,\md\x$ and $t\mapsto\rb_i(t)$, in $\R_+$ (up to $t=0$), for $i\in\{1,2\}$. Now, for any $i\neq j\in\{1,2\}$, $0<\epsilon<t$ and $R>0$, integrating~\eqref{eq:sys_main1} over $(\epsilon,t)\times B(0,R)$ yields:
$$\baa{rcl}
\displaystyle\int_{B(0,R)}u_i(t,\x)\,\md\x-\int_{B(0,R)}u_i(\epsilon,\x)\,\md\x & = & \displaystyle\frac{\mu^2}{2}\int_\epsilon^t\int_{\partial B(0,R)}\nu\cdot\nabla u_i(s,\x)\,\md\sigma(\x)\,\md s\\
& & \displaystyle+\int_\epsilon^t\int_{B(0,R)}r_i(\x)\,u_i(s,\x)\,\md\x\,\md s\\
& & \displaystyle+\int_\epsilon^t\int_{B(0,R)}\big(\delta_{i,j}u_j(s,\x)-\delta_{i,i}u_i(s,\x)\big)\,\md\x\,\md s,\eaa$$
where $\nu$ and $\md\sigma(\x)$ denote the outward normal and surface measure on $\partial B(0,R)$. From~\eqref{eq:sys_main1},~\eqref{eq:maj1:dens} and~\eqref{eq:ineghi}, together with~(H3) and standard parabolic estimates, it follows that $\|\x\|^{n+1}u_i(s,\x)\to0$ and $\|\x\|^{n-1}\|\nabla u_i(s,\x)\|\to0$ as $\|\x\|\to+\infty$, uniformly for $s\in[\epsilon,t]$. Therefore, by passing to the limit $R\to+\infty$ in the above displayed equality, one gets that:
\[N_i(t)-N_i(\epsilon)=\int_\epsilon^t\rb_i(s)\,N_i(s)\,\md s+\int_\epsilon^t\big(\delta_{i,j}N_j(s)-\delta_{i,i}N_i(s)\big)\,\md s,\]
where we also used Lebesgue's dominated convergence theorem, formula~\eqref{def:m(x)} and the continuity of the map $s\mapsto\int_{\R^n}m_i(\x)\,u_i(s,\x)\,\md\x$ in $\R_+$. Using the continuity of $N_i$, $N_j$ and $\rb_i$ in $\R_+$, the passage to the limit $\epsilon\to0^+$ yields:
\[N_i(t)-N_i(0)=\int_0^t\rb_i(s)\,N_i(s)\,\md s+\int_0^t\big(\delta_{i,j}N_j(s)-\delta_{i,i}N_i(s)\big)\,\md s.\]
Hence, each function $N_i$ is of class $C^1(\R_+)$ and the pair $(N_1,N_2)$ satisfies~\eqref{eq:ODEmalthusbis}.
\vskip 0.3cm
\noindent{}Part~(ii). We now show the symmetry property of the solutions of~\eqref{eq:sys_main1} with fitnesses~\eqref{def:ri}, still for the first type~\eqref{eq:f_type1}, with~$(\dens_1,\dens_2)$ given as above in part~(i) and initial conditions now satisfying~\eqref{hyp:HS} as well. With these assumptions, it follows that the pair of functions $(U_1,U_2)$ defined by:
\begin{equation*}
\forall\,t\in\R_+, \ \forall\,\x \in \R^n, \ \  
(U_1(t,\x),U_2(t,\x))= (\dens_2(t,\iota(\x)),\dens_1(t,\iota(\x))),
\end{equation*}
with $\iota$ as in \eqref{dfn:iota}, is a $C^{1,2}([0,+\infty)\times\R^n)^2$ solution of the Cauchy problem~\eqref{eq:sys_main1}. Furthermore, each component $U_i$ is positive in $(0,+\infty)\times\R^n$, bounded in $(0,T)\times\R^n$ for every $T>0$, and converges to $0$ as $\|\x\|\to+\infty$ locally uniformly in $t\in\R_+$. By uniqueness of such solutions and by~\eqref{hyp:HS}, one gets that $U_1(t,\x) =\dens_1(t,\x)$ and $U_2(t,\x) = \dens_2(t,\x)$ for all $(t,\x)\in \R_+\times\R ^n$, and so~$\dens_1(t,\x) = \dens_2(t,\iota(\x)).$ The equation~\eqref{pde:denssym} then readily follows from this equality. Moreover the population sizes at time $t\ge 0$ satisfy:
\[
N_1(t) = \int_{\R^n} \dens_1(t,\x) \,\md\x = \int_{\R^n} \dens_2(t,\x)\,\md\x=N_2(t),
\]
and the mean fitnesses are also such that $\rb_1(t)=\rb_2(t)$ for all $t\ge 0$.

In order to complete the proof of Theorem~\ref{thm:existunidens}, we now assume that the fitnesses are given by~\eqref{def:ri} and we derive an equivalence between the problem~\eqref{eq:sys_main1} in the symmetric case~\eqref{hyp:HS} with $f_1,f_2$ of the first type~\eqref{eq:f_type1}, and the problem~\eqref{eq:sys_main1} with~$f_1,f_2$ of the second type~\eqref{eq:f_type2}, still in the symmetric case~\eqref{hyp:HS}. Firstly, assume that $f_1,f_2$ are of the first type~\eqref{eq:f_type1}, and let $\dens_i$,  $N_i$ and $\rb_i$ be defined by the first part of the present proof, for $i\in\{1,2\}$. From~\eqref{hyp:HS} and the previous paragraph, we know that $\rb_1(t)= \rb_2(t)=:\rb(t)$ and $N_1(t)= N_2(t)=:N(t)>0$, with $N'(t)=\rb(t)\,N(t)$, for all $t\ge 0$. Let $\tN(t)$ be the solution of the ODE:
$$\tN'(t)=\rb(t)\, \tN(t) - \tN(t)^2,$$
with $\tN(0)=N(0)>0$. Since $\rb$ is continuous in $\R_+$, the function $\tN$ is well defined, positive, and of class~$C^1$ in $\R_+$. Define, for $i\in\{1,2\}$, the functions:
$$\forall\,t\in\R_+, \ \forall\,\x \in \R^n,\ \  \tu_i(t,\x)=\frac{\tN(t)}{N(t)}\, \dens_i(t,\x),$$
where the functions $u_i$ are recalled to satisfy~\eqref{eq:sys_main1} with $f_1,f_2$ of the first type~\eqref{eq:f_type1}. The pair $(\tu_1,\tu_2)$ is of class $C^{1,2}(\R_+\times\R^n)^2$, it is locally bounded in time, it converges to $(0,0)$ as $\|\x\|\to+\infty$ locally uniformly in $t\in\R_+$, and it has the same initial condition as the pair $(\dens_1,\dens_2)$. Moreover, for all~$t\ge 0$ and $\x\in \R^n$, we have:
$$\frac{\tN(t)}{N(t)}\,\partial_t\dens_i(t,\x)={\dd^2\over 2}\Delta \tu_i(t,\x)+r_i(\x)\tu_i(t,\x) + \delta\, [\tu_j(t,\x)-\tu_i(t,\x)],\footnote{We use in the last term the fact that the proportionality factor between $\tu_i$ and $\dens_i$ is the same for $i\in\{1,2\}$.}$$
and:
\begin{align*}
\partial_t \tu_i(t,\x)& =\frac{\tN(t)}{N(t)}\,\partial_t\dens_i(t,\x)+\lp\frac{\tN'(t)\, N(t)-\tN(t) \, N'(t)}{N^2(t)}\rp \dens_i(t,\x),\\
& =\frac{\tN(t)}{N(t)}\,\partial_t\dens_i(t,\x)-\frac{\tN(t)^2}{N(t)}\,\dens_i(t,\x)=\frac{\tN(t)}{N(t)}\,\partial_t\dens_i(t,\x)-\tN(t)\tu_i(t,\x).
\end{align*}
The functions $\tu_i$ thus satisfy (with $i,j\in\{1,2\}$ and $i\neq j$):
$$\partial_t\tu_i(t,\x)={\dd^2\over 2}\Delta \tu_i(t,\x)+ \left[ r_i(\x)-\tN(t)\right] \, \tu_i(t,\x)+ \delta [\tu_j(t,\x)-\tu_i(t,\x)],$$
for all $t\ge0$ and $\x\in\R^n$, and, as: $$\int_{\R^n}\tu_i(t,\x)\, \md\x=\tN(t),$$
for all $t\ge0$ and $i\in\{1,2\}$, the functions $\tu_i$ then solve~\eqref{eq:sys_main1}, with $f_1,f_2$ of the second type~\eqref{eq:f_type2}. These solutions $\tu_i$ are also symmetric, in the sense that $\tu_1(t,\x)=\tu_2(t,\iota(\x))$ for all $t\ge0$ and $\x\in\R^n$, and they are positive in $(0,+\infty)\times\R^n$. Notice finally that:
$$\widetilde{r}(t):=\frac{1}{\tN(t)}\int_{\R^n}r_i(\x)\,\tu_i(t,\x)\,\md\x=\rb(t),$$
for all $t\ge0$ and $i\in\{1,2\}$.

Conversely, assume that $(\tu_1,\tu_2)$ is a symmetric $C^{1,2}(\R_+\times\R^n)^2$ locally bounded in time solution of~\eqref{eq:sys_main1} and converging to $(0,0)$ as $\|\x\|\to+\infty$ locally uniformly in $t\in\R_+$, with $f_1,f_2$ of the second type~\eqref{eq:f_type2} and with a continuous associated population size $\tN(t)$ in each habitat, such that $\tN(0)>0$. Since the system satisfied by $(\tu_1,\tu_2)$ can also be viewed as a linear cooperative system (with additional diagonal term $-\tN_i(t)\,\tu_i(t,x)$), the weak and strong comparison principle applied with respect to the trivial solution $(0,0)$ imply that the functions $\tu_i$ are then positive in~$(0,+\infty)\times\R^n$. Therefore, the population size $\tN(t)$ is positive and $f_i(\x,[\tu_i])\Bk\le r_i(\x)\,\tu_i(t,\x)$ for all $t\ge0$ and $\x\in\R^n$. As a consequence, the pair $(\tu_1,\tu_2)$ is then a subsolution of the cooperative system~\eqref{eq:sys_main1} with growth functions of the first type~\eqref{eq:f_type1}. Since the maximum principle holds for the latter system, one infers that the functions $\tu_i$ satisfy similar bounds as~\eqref{eq:maj1:dens} and~\eqref{eq:ineghi} above for the solutions $u_i$ in the first type~\eqref{eq:f_type1}. By arguing as above, it follows that the mean fitness $t\mapsto\widetilde{r}(t)=\tN(t)^{-1}\int_{\R^n}r_i(\x)\,\tu_i(t,\x)\,\md\x$ is continuous in $\R_+$ and independent of $i\in\{1,2\}$, and that population size $\tN$ is of class $C^1(\R_+)$ and satisfies~\eqref{eq:ODElogistic} (due to the additional term $-\tN(t)\,\tu_i(t,\x)$ in the right-hand side of the equation satisfied by $\tu_i$). Finally, by inverting all the calculations of the previous paragraph and by defining $N(t)$ as the solution of $N'(t)=\widetilde{r}(t)\,N(t)$ with $N(0)=\tN(0)$, one gets that the pair $(\dens_1,\dens_2)$ defined by:
$$\forall\,t\in\R_+,\ \forall\,\x \in \R^n,\ \  \dens_i(t,\x)=\frac{N(t)}{\tN(t)}\, \tu_i(t,\x),$$
is a symmetric solution of~\eqref{eq:sys_main1} satisfying the conditions of Theorem~\ref{thm:existunidens} with growth functions $f_1,f_2$ of the first type~\eqref{eq:f_type1}. The uniqueness result for the solutions in the first type~\eqref{eq:f_type1} then leads to the uniqueness of the symmetric solutions of~\eqref{eq:sys_main1} for growth functions of the second type~\eqref{eq:f_type2}. The proof of Theorem~\ref{thm:existunidens} is thereby complete.
\qed

\subsection{Large time behaviour}\label{sec:32}

This section is devoted to the proofs of Theorems~\ref{th:extinction_vs_explo:lambda} and~\ref{th:extinction_vs_explo:logistic}. Before that, we state an auxiliary lemma on the existence of positive eigenfunctions of the operator $\A$ defined in~\eqref{eq:defAbis}, associated to problem~\eqref{eq:sys_main1bis} with fitnesses~\eqref{def:ribis}.

\begin{lemma}\label{lem:fpp}
There exists a pair of positive eigenfunctions $(\varphi_1, \varphi_2)\in\big(C^\infty_0(\R^n)\cap L^1(\R^n)\big)^2$ such that $\A\, (\varphi_1,\varphi_2)= \lambda (\varphi_1,\varphi_2)$ in $\R^n$, with $\A$ defined by~\eqref{eq:defAbis} with fitnesses~\eqref{def:ribis}, and with $\lambda$ defined by~\eqref{eq:def_lambda}. Furthermore, this pair $(\varphi_1,\varphi_2)$ is unique up to multiplication by a positive constant. Lastly, for the particular case of the operator $\A$ defined in~\eqref{eq:defA} with fitnesses~\eqref{def:ri}, the functions~$\varphi_1$ and~$\varphi_2$ satisfy $\varphi_1=\varphi_2\circ\iota$ in $\R^n$.
\end{lemma}

The proof of Lemma~\ref{lem:fpp} is postponed after that of Theorem~\ref{th:extinction_vs_explo:lambda}.

\vskip 0.3cm

\noindent \textit{Proof of Theorem~\ref{th:extinction_vs_explo:lambda}}. Let $\densbf=(\dens_1,\dens_2)$ be the solution of~\eqref{eq:sys_main1bis} and~\eqref{def:ribis} given by Theorem~\ref{thm:existunidens}-(i), with an initial condition $\densbf^0=(u^0_1,u^0_2)$ satisfying~(H1)-(H3), for $f_1,f_2$ of the first type~\eqref{eq:f_type1}. Let~$N_1(t)$ and $N_2(t)$ be its population sizes in each habitat, at time $t\ge0$. For $R>0$, let $(\varphi_1^{R}, \varphi_2^{R})\in C^\infty_0(\overline{B(0,R)})^2$  and $\lambda^{R}$ be the principal eigenfunctions and eigenvalue of the operator~$\A$ defined by \eqref{eq:defAbis}. Finally, let $\lambda$ be given by~\eqref{eq:def_lambda}. We consider the cases $\lambda<0$ and $\lambda\ge0$ separately.

{\it First case: assume that $\lambda<0$.} From assumptions~(H2)-(H3), we know that $u^0_i\ge0$ and~$u^0_i \not \equiv 0$ in $\R^n$, for $i=1,2$, and, from Theorem~\ref{thm:existunidens}-(i), $\dens_i(1,\cdot)>0$ in $\R^n$ for each $i\in\{1,2\}$. As $\lim_{R\to+\infty}\lambda^{R}=\lambda<0$, we can fix $R>0$ such that $\lambda^{R}<0$. Let $\underline K>0$ be such that $\underline K\,e^{-\lambda^{R}}(\varphi_1^{R}, \varphi_2^{R})\le \densbf(1,\cdot)$ in~$\overline{B(0,R)}$. Set~$\underline H(t,\x)=(\underline H_1,\underline H_2)(t,\x):=\underline K\, e^{-\lambda^{R}t}(\varphi_1^{R}(\x), \varphi_2^{R}(\x))$ for $t\ge 1$ and $\x\in\overline{B(0,R)}$. In particular, $\underline H(1,\cdot)\le\densbf(1,\cdot)$ in $\overline{B(0,R)}$. We have, for all $t\ge 1,$ and $i\neq j\in\{1,2\}$:
\be\label{eq:Hi}
\left\{\baa{ll}
\partial_t\underline H_i=\displaystyle\frac{\mu^2}{2} \Delta\underline H_i +r_i(\x)\,\underline H_i +\delta_{i,j}\underline H_j-\delta_{i,i}\underline H_i, & \hbox{ in }\overline{B(0,R)},\vspace{3pt}\\
(\underline H_1,\underline H_2)(t,\cdot)=(0,0), & \hbox{ on }\partial B(0,R).\eaa\right.
\ee
As the pair $(\dens_1(t,\cdot),\dens_2(t,\cdot))$ satisfies the same equation in $\overline{B(0,R)}$ and is positive in $\R^n$ for each $t\ge 1$ and therefore on $\partial B(0,R)$, the maximum principle applied to this cooperative system implies that $\dens_i(t,\x)\ge \underline H_i(t,\x)$ for all $t\ge 1$, $\x\in\overline{B(0,R)}$ and $i\in\{1,2\}$. Integrating over $B(0,R)$ the above inequality and using the positivity of $u_i$, we get:
$$N_i(t)\ge \underline K\, e^{-\lambda^{R}t}\,\|\varphi_i^{R}\|_{L^1(B(0,R))},\ \hbox{ for all }t\ge 1\hbox{ and }i\in\{1,2\}.$$
Since $\lambda^{R}<0$, this implies that $N_i(t)\to +\infty$ as $t\to +\infty$ and this shows part~(i) of Theorem~\ref{th:extinction_vs_explo:lambda}. Notice that the above inequality also implies that:
\be\label{ineqNi}
\liminf_{t\to+\infty}\frac{\ln N_i(t)}{t}>0.
\ee

{\it Second case: assume that $\lambda\ge0$.} Assume also that the initial conditions $u^0_1$ and $u^0_2$ are compactly supported. Then, there is $\overline K>0$ large enough, one has $\overline K\,(\varphi_1, \varphi_2) \ge \densbf(0,\cdot)$ in $\R^n$. Set $\overline H(t,\x)=(\overline H_1,\overline H_2)(t,\x):= \overline K\, e^{-\lambda t}(\varphi_1(\x), \varphi_2(\x))$ for $t\ge 0$ and $\x \in \R^n$. As for~\eqref{eq:Hi}, the function~$\overline H$ satisfies the same cooperative system~\eqref{eq:sys_main1bis} as $\densbf$ in $\R_+\times\R^n$, but with a larger initial condition. The comparison principle thus implies that, for $i\in\{1,2\}$ and $t\ge 0$:
\begin{equation} \label{eq:ineg_uihi}
0\le\dens_i(t,\x)\le \overline H_i(t,\x)=\overline K\, e^{-\lambda t}\varphi_i(\x),\ \hbox{ for all }\x \in \R^n.
\end{equation}
As the functions $\varphi_i$ belong to $L^1(\R^n)$, integrating \eqref{eq:ineg_uihi} over $\R^n$ yields $\lim_{t\to+\infty}N_i(t)=0$ if~$\lambda>0$. If~$\lambda=0$,~\eqref{eq:ineg_uihi} implies that:
\[\displaystyle\limsup_{t\to+\infty}N_i(t)<+\infty,\ \hbox{ for }i=1,2.
\]
Furthermore, in that case, for every~$C>0$, $C\, (\varphi_1, \varphi_2)$ is a pair of positive stationary solutions of~\eqref{eq:sys_main1bis} with finite population sizes. That shows parts~(ii) and~(iii) of Theorem~\ref{th:extinction_vs_explo:lambda} and the proof of Theorem~\ref{th:extinction_vs_explo:lambda} is thereby complete.
\qed

\vskip 0.3cm

\noindent \textit{Proof of Lemma~\ref{lem:fpp}.} For $R>0$, the functions $(\varphi_1^{R}, \varphi_2^{R})\in \big(C^\infty_0(\overline{B(0,R)})\big)^2$ satisfy:
$$\frac{\mu^2}{2}\Delta \varphi_i^{R}+(\lambda^{R} -\delta_{i,i}+r_i)\,\varphi_i^{R} + \delta_{i,j}\varphi_j^{R}=0\ \hbox{ in }\overline{B(0,R)},$$
with $i\neq j\in\{1,2\}$. As the eigenvalues $\lambda^{R}$ are nonincreasing with respect to $R$ and not smaller than $\min(-r_{\max,1}+\delta_{1,1}-\delta_{1,2},-r_{\max,2}+\delta_{2,2}-\delta_{2,1})$ by~\eqref{ineqlambda}, we have:
\[\min(-r_{\max,1}+\delta_{1,1}-\delta_{1,2},-r_{\max,2}+\delta_{2,2}-\delta_{2,1})\le\lambda^{R}\le\lambda^{2}\]
for all $R\ge2$. For every $R'\ge1$, it then follows from the Harnack inequality in~\cite[Theorem~2]{Sir09} (applied here with $\Omega=B(0,2R')$) that there is a positive constant $C(R')$ such that:
$$\max_{\x \in\overline{B(0,R')},\,i\in\{1,2\}}\,\varphi^{R}_i(\x)\le C(R') \,\min_{\x\in\overline{B(0,R')},\,i\in\{1,2\}}\varphi^{R}_i(\x),$$
for all $R\ge 2R'$. Without loss of generality, up to multiplication by a positive constant, we assume the normalization condition:
\[\min\big(\varphi_1^{R}(0),\varphi_2^R(0)\big)=1.\]
Thus, we get:
$$0< \varphi^{R}_i(\x) \le C(R'),\ \hbox{ for all }\x\in\overline{B(0,R')},\ i\in\{1,2\},\hbox{ and }R\ge 2R'\ge2.$$
Standard elliptic estimates then imply that, for every $\theta\in[0,1)$, and for every $R'\ge1$, the functions~$\varphi^{R}_i$ are bounded in $C^{2,\theta}(\overline{B(0,R')})$, independently of $R\in[2R',+\infty)$. Thus Sobolev's injections imply that, up to the extraction of a subsequence, $\varphi^{R}_i\to \varphi_i$ in $C^2_{loc}(\R^n)$ as $R\to+\infty$, where the functions $\varphi_i$ satisfy
$\A\,(\varphi_1,\varphi_2)= \lambda (\varphi_1,\varphi_2)$,  are nonnegative and such that $\varphi_i(0)=1$ for $i\in\{1,2\}$. From the (scalar) strong elliptic maximum principle, the functions $\varphi_i$ are positive in $\R^n$. Furthermore, again from standard elliptic estimates, they are of class $C^\infty(\R^n)$. Notice also from Remark~\ref{rem:phi_symm} that, for the operator $\A$ defined in~\eqref{eq:defA} associated to problem~\eqref{eq:sys_main1} with fitnesses~\eqref{def:ri}, the functions $\varphi^R_i$ satisfy $\varphi^R_1=\varphi^R_2\circ\iota$ in $\overline{B(0,R)}$, hence $\varphi_1=\varphi_2\circ\iota$ in $\R^n$ in that case.

To show that the eigenfunctions $\varphi_i$ are in $L^1(\R^n)$ and converge to $0$ at infinity, we use the fact that the potentials $r_i$ in~\eqref{def:ribis} are confining. In particular, we fix $R'_0\ge1$ large enough such that, all $R\ge2R'_0$, there holds:
$$\max\big(r_1(\x)+\delta_{2,1}-\delta_{1,1}+\lambda^{R},r_2(\x)+\delta_{1,2}-\delta_{2,2}+\lambda^{R}\big)<-\frac{\|\x\|^2}{4},$$
for all $\x\in\overline{B(0,R)}\setminus B(0,R'_0)$ and $i\in\{1,2\}$, hence:
$$-\frac{\mu^2}{2}\Delta (\varphi_1^{R}+\varphi_2^{R})(\x) + \frac{\|\x\|^2}{4} \, (\varphi_1^{R}+\varphi_2^{R})(\x)< 0\ 
\hbox{ in }\overline{B(0,R)} \setminus B(0,R'_0).$$
For any such $R$, since $\max_{\partial B(0,R'_0)}\varphi_i^{R}\le C(R'_0)$ and $\varphi_1^{R}+\varphi_2^{R}=0$ on $\partial B(0,R)$, the maximum principle implies that  $\varphi_1^{R}+\varphi_2^{R}\le w$ in $\overline{B(0,R)} \setminus B(0,R'_0)$, where $w$ denotes the solution of the equation $-(\mu^2/2)\Delta w(\x)+(\|\x\|^2/4) w(x)=0$ in $\overline{B(0,R)}\setminus B(0,R'_0)$  with the boundary conditions $w=0$ on~$\partial B(0,R)$ and $w=2 \, C(R'_0)$ on~$\partial B(0,R'_0)$. Consequently,
$$\varphi_1^{R}(x)+\varphi_2^{R}(x)\le w(x)\le 2\,C(R'_0)\,e^{({R'_0}^2-\|\x\|^2)/\sqrt{8\mu^2}},\quad \hbox{ for all }\x\in\overline{B(0,R)}\setminus B(0,R'_0),$$
and for all $R\ge 2R'_0$. Thus, the same inequality holds for the functions  $\varphi_1+\varphi_2$ in $\R^n\setminus B(0,R'_0)$. This implies in particular that the eigenfunctions $\varphi_i$ belong to $L^1(\R^n)$ and converge to $0$ at infinity.

Lastly, since for any $\lambda\in\R$ the weak maximum principle holds outside a large ball for the system $\A(\phi_1,\phi_2)=\lambda(\phi_1,\phi_2)$ in the class of $C^2_0(\R^n)$ functions (namely, there is $\rho>0$ such that, if $\phi_1,\phi_2\in C^2_0(\R^n)$ satisfy $\A(\phi_1,\phi_2)\le\lambda(\phi_1,\phi_2)$ in $\R^n\setminus B(0,\rho)$ and $(\phi_1,\phi_2)\le(0,0)$ on $\partial B(0,\rho)$, then $(\phi_1,\phi_2)\le(0,0)$ in $\R^n\setminus B(0,\rho)$) and since the strong maximum principle holds as well in any connected open subset $\Omega\subset\R^n$ (namely, if $\phi_1,\phi_2\in C^2(\Omega)$ are such that $\A(\phi_1,\phi_2)\le\lambda(\phi_1,\phi_2)$ and $(\phi_1,\phi_2)\le(0,0)$ in $\Omega$ with $\phi_i(x_0)=0$ for some $i\in\{1,2\}$ and $x_0\in\Omega$, then $(\phi_1,\phi_2)\equiv(0,0)$ in~$\Omega$), it follows with similar arguments as in~\cite{BNV} that the pair of eigenfunctions $(\varphi_1,\varphi_2)$ constructed above is unique, up to multiplication by a positive constant, in the class of $C^2_0(\R^n)$ eigenfunctions. Moreover, the eigenvalue $\lambda$ is the unique eigenvalue associated with a pair of positive eigenfunctions. The proof of Lemma~\ref{lem:fpp} is thereby complete.
\qed

\vskip 0.3cm

\noindent \textit{Proof of Theorem~\ref{th:extinction_vs_explo:logistic}}.
Let $\densbf=(\dens_1,\dens_2)$ be the unique symmetric solution of~\eqref{eq:sys_main1} and~\eqref{def:ri} given by Theo\-rem~\ref{thm:existunidens}-(ii), for $f_1,f_2$ of the second type~\eqref{eq:f_type2} and for an initial condition $\densbf^0=(u^0,u^0\circ\iota)$ satisfying~\eqref{hyp:HS} and~(H1)-(H3). Let $N(t):=N_1(t)=N_2(t)$ be its population size given by~\eqref{dfn:pop_size_pde} and $\rb(t):=\rb_1(t)=\rb_2(t)$ be its mean fitness given by~\eqref{dfn:meanfit}, at each time $t\ge0$. From Theorem~\ref{thm:existunidens}, the densities $u_i$ are positive in $(0,+\infty)\times\R^n$, the function $\rb$ is continuous in $\R_+$, the function $N$ is positive and of class $C^1$ in $\R_+$, and $N'(t)=\rb(t)N(t)-N(t)^2$ for all $t\in\R_+$.

Let also $\widetilde{\densbf}=(\tu_1,\tu_2)$ be the unique, symmetric, solution of~\eqref{eq:sys_main1} and~\eqref{def:ri} given by Theorem~\ref{thm:existunidens}, for~$f_1,f_2$ of the first type~\eqref{eq:f_type1}, with the same initial condition $\densbf^0=(u^0,u^0\circ\iota)$ as $\densbf$. Let $\tN(t):=\tN_1(t)=\tN_2(t)$ be its population size and $\widetilde{r}(t):=\widetilde{r}_1(t)=\widetilde{r}_2(t)$ be its mean fitness, at each time $t\ge0$. From Theorem~\ref{thm:existunidens}, the densities $\tu_i$ are positive in $(0,+\infty)\times\R^n$, the function $\widetilde{r}$ is continuous in $\R_+$, the function $\tN$ is positive and of class $C^1$ in $\R_+$, and $\tN'(t)=\widetilde{r}(t)\tN(t)$ for all $t\in\R_+$. 

The correspondence between the symmetric solutions of~\eqref{eq:sys_main1} for both types~\eqref{eq:f_type1} and~\eqref{eq:f_type2}, shown in the proof of part~(ii) of Theorem~\ref{thm:existunidens}, implies that:
$$\widetilde{\densbf}(t,\x)=\frac{\tN(t)}{N(t)}\,\densbf(t,\x),\ \hbox{ for all }t\ge0\hbox{ and }\x\in\R^n,$$
hence $\widetilde{r}(t)=\rb(t)$ for all $t\ge0$. Therefore, we have:
\be\label{N'}
\frac{\tN'(t)}{\tN(t)}=\frac{N'(t)}{N(t)}+N(t),
\ee
for all $t\ge0$. Integrating this equality and using $\tN(0)=N(0)$ yields:
\be\label{eq:NtN}
N(t)=\frac{\tN(t)}{\displaystyle 1+\int_0^t\tN(s)\,\md s},\quad \hbox{ for all }t\ge0.
\ee

Let now $(\varphi_1, \varphi_2)$ be defined by Lemma~\ref{lem:fpp} with the normalization $\|\varphi_i\|_{L^1(\R^n)} =1$ (we recall that the functions $\varphi_1$ and $\varphi_2$ are here such that $\varphi_1=\varphi_2\circ\iota$ in $\R^n$). Set:
$$H(t,\x)=(H_1(t,\x),H_2(t,\x)):= e^{-\lambda t}(\varphi_1(\x), \varphi_2(\x)),$$
for $t\ge0$ and $\x\in\R^n$. As in the proof of Theorem~\ref{th:extinction_vs_explo:lambda}, the function $H$ satisfies~\eqref{eq:sys_main1} with growth functions $f_i$ of the first type~\eqref{eq:f_type1}. We then treat separately the cases $\lambda\ge0$ and $\lambda<0$.

{\it First case: Assume that $\lambda\ge0$.} Assume also in this case that $u^0$ is compactly supported. Then there is $K>0$ such that $\widetilde{\densbf}(0,\cdot)=\densbf(0,\cdot)=\densbf^0\le K\,H(0,\cdot)$ in $\R^n$ and the maximum principle applied to the cooperative system~\eqref{eq:sys_main1} with growth functions of the first type~\eqref{eq:f_type1} implies that~$\widetilde{\densbf}(t,\cdot)\le K\,H(t,\cdot)$ in $\R^n$ for all $t\ge0$, hence $\tN(t)\le K\,e^{-\lambda t}$ for all $t\ge0$. From~\eqref{eq:NtN} and the positivity of $N$ and $\tN$, one immediately infers that $N(t)\to0$ as $t\to+\infty$ if $\lambda>0$.

Consider now the sub-case $\lambda=0$. The previous observations imply that $\tN$ is bounded in~$\R_+$. Furthermore, on the one hand, if the integral $\int_0^{+\infty}\tN(s)\,\md s$ diverges, then formula~\eqref{eq:NtN} and the boundedness of $\tN$ imply that $N(t)\to0$ as $t\to+\infty$. On the other hand, if the integral $\int_0^{+\infty}\tN(s)\,\md s$ converges, then the boundedness of the function $\tN'=\widetilde{r}\,\tN$ in $\R_+$ (which itself follows from the inequalities $0\le\widetilde{\densbf}(t,\cdot)\le K\,H(t,\cdot)=K\,(\varphi_1,\varphi_2)$ in $\R^n$ and the exponential decay at infinity of the eigenfunctions~$\varphi_i$ given the proof of Lemma~\ref{lem:fpp}) implies that $\tN(t)\to0$ as~$t\to+\infty$, and finally $N(t)\to0$ as~$t\to+\infty$ by~\eqref{eq:NtN}.

{\it Second case: Assume that $\lambda<0$.} Assume also in this case that $\densbf^0$ is trapped between two positive multiples of the eigenfunctions $(\varphi_1,\varphi_2)$, namely, there exist $0<K_1\le K_2$ such that:
\be\label{hypu0}
K_1\,(\varphi_1,\varphi_2)\le\densbf^0\le K_2\,(\varphi_1,\varphi_2)\ \hbox{ in $\R^n$}.
\ee
Thus, $K_1\,H(0,\cdot)\le\widetilde{\densbf}(0,\cdot)=\densbf^0\le K_2\,H(0,\cdot)$ in $\R^n$ and the maximum principle applied to the coope\-ra\-tive system~\eqref{eq:sys_main1} with growth functions of the first type~\eqref{eq:f_type1} implies that:
\[
K_1\,H(t,\cdot)\le\widetilde{\densbf}(t,\cdot)\le K_2\,H(t,\cdot)\ \text{ in } \R^n \text{, for all } t\ge0.
\]In particular, $K_1\,e^{-\lambda t}\le\tN(t)\le K_2\,e^{-\lambda t}$ for all $t\ge0$. Together with~\eqref{eq:NtN} and the negativity of~$\lambda$, one concludes that:
$$0<\frac{K_1}{K_2}\,|\lambda|\le\liminf_{t\to+\infty}N(t)\le\limsup_{t\to+\infty}N(t)\le\frac{K_2}{K_1}\,|\lambda|<+\infty.$$
Furthermore, for initial conditions that satisfy~\eqref{hyp:HS} and~(H1)-(H3) but may not satisfy~\eqref{hypu0}, one knows from the proof of Theorem~\ref{th:extinction_vs_explo:lambda}, namely from~\eqref{ineqNi}, that $\liminf_{t\to+\infty}(\ln\tN(t))/t>0$.\footnote{Notice that in formula~\eqref{ineqNi} of the proof of Theorem~\ref{th:extinction_vs_explo:lambda}, $N(t)$ was the population size of the solution for growth functions of the first type~\eqref{eq:f_type1}, whereas here this population size is called $\tN(t)$.} On the other hand, integrating~\eqref{N'} over $(0,t)$ and using $\tN(0)=N(0)$ leads to:
$$\ln\tN(t)=\ln N(t)+\int_0^tN(s)\,\md s,$$
for every $t>0$. Hence, $\liminf_{t\to+\infty}N(t)>0$, since otherwise the right-hand side of the above equality would be not larger than $o(t)$ as $t\to+\infty$, then contradicting $\liminf_{t\to+\infty}(\ln\tN(t))/t>0$. The proof of Theorem~\ref{th:extinction_vs_explo:logistic} is thereby complete.
\qed

\subsection{Dependence with respect to the parameters}\label{sec:33}

\noindent \textit{Proof of Proposition~\ref{prop:continuity_lambda}}.
{\it Step 1: general properties of $\lambda_{\delta,m_D,\mu,\rmax}$.} We recall that $\lambda=\lambda_{\delta,m_D,\mu,\rmax}$ denotes the principal eigenvalue, given in~\eqref{eq:def_lambda}, for the operator $\A$ defined by~\eqref{eq:defA} in $\R^n$, with fitnesses~\eqref{def:ri}. We also recall that $m_D=2\beta^2$, with~$\beta$ given in~\eqref{def:opt_beta}. In this Step~1, $\beta$ can actually be any real number. Using the confining properties of the fitnesses~$r_i(\x)$, it follows from Lemma~\ref{lem:fpp} and elementary arguments that, for every $\delta>0$, $m_D\ge0$, $\mu>0$ and~$\rmax\in\R$,
\be\label{lambdadelta}
\lambda_{\delta,m_D,\mu,\rmax}=\ \min_{\varphi\in\mathcal{E}}\,\mathcal{R}_{\delta,\beta,\mu,\rmax}(\varphi),
\ee
where:
$$\mathcal{E}=\big\{\varphi\in H^1(\R^n),\,\x\mapsto\|\x\|\,\varphi(\x)\,\in L^2(\R^n),\ \|\varphi\|_{L^2(\R^n)}=1\big\},$$
and:
$$\mathcal{R}_{\delta,\beta,\mu,\rmax}(\varphi)=\frac{\mu^2}{2}\int_{\R^n}\|\nabla\varphi(\x)\|^2\,\md\x-\int_{\R^n}r_1(\x)\,\varphi(\x)^2\,\md\x+\delta\int_{\R^n}\big(\varphi(\x)^2-\varphi(\x)\,\varphi(\iota(\x))\big)\,\md\x,$$
and the minimum of $\mathcal{R}_{\delta,\beta,\mu,\rmax}$ in~\eqref{lambdadelta} is reached only by $\pm\varphi_1$, where $\varphi_1$ is the positive eigenfunction given in Lemma~\ref{lem:fpp}, normalized so that $\|\varphi_1\|_{L^2(\R^n)}=1$ (notice that this function $\varphi_1$ belongs to $\mathcal{E}$ from the bounds derived in the proof of Lemma~\ref{lem:fpp}). One can also write:
$$\mathcal{R}_{\delta,\beta,\mu,\rmax}(\varphi)=-\rmax+\mathcal{S}_{\delta,\beta,\mu}(\varphi)=-\rmax+\frac{\beta^2}{2}+\mathcal{T}_{\delta,\beta,\mu}(\varphi),$$
with:
\begin{multline*}
\mathcal{S}_{\delta,\beta,\mu}(\varphi)  =  \displaystyle\frac{\mu^2}{2}\int_{\R^n}\|\nabla\varphi(\x)\|^2\,\md\x+\int_{\R^n}\frac{\|(x_1+\beta,x_2,\cdots,x_n)\|^2}{2}\,\varphi(\x)^2\,\md\x\vspace{3pt}\\
\displaystyle+\delta\int_{\R^n}\big(\varphi(\x)^2-\varphi(\x)\,\varphi(\iota(\x))\big)\,\md\x\end{multline*}
and:
$$\mathcal{T}_{\delta,\beta,\mu}(\varphi)=\frac{\mu^2}{2}\!\int_{\R^n}\!\!\|\nabla\varphi(\x)\|^2\,\md\x+\int_{\R^n}\!\!\Big(\beta x_1+\frac{\|\x\|^2}{2}\Big)\,\varphi(\x)^2\,\md\x+\delta\!\int_{\R^n}\!\!\big(\varphi(\x)^2-\varphi(\x)\,\varphi(\iota(\x))\big)\,\md\x.$$
Therefore,
\be\label{lambdaT}
\lambda_{\delta,m_D,\mu,\rmax}=-\rmax+\frac{\beta^2}{2}+\min_{\varphi\in\mathcal{E}}\,\mathcal{T}_{\delta,\beta,\mu}(\varphi).
\ee
Notice immediately that:
$$\lambda_{\delta,m_D,\mu,\rmax}=-\rmax+\lambda_{\delta,m_D,\mu,0},$$
for all $(\delta,m_D,\mu,\rmax)\in(0,+\infty)\times[0,+\infty)\times(0,+\infty)\times\R$. Furthermore, since $\mathcal{T}_{\delta,\beta,\mu}(\varphi)$ is affine with respect to $(\delta,\beta,\mu^2)\in(0,+\infty)\times\R\times(0,+\infty)$ for each $\varphi\in\mathcal{E}$, one also infers that the function $(\delta,\beta,\mu^2)\mapsto\min_{\varphi\in\mathcal{E}}\,\mathcal{T}_{\delta,\beta,\mu}(\varphi)$ is concave in $(0,+\infty)\times\R\times(0,+\infty)$, hence it is continuous in this set. This together with~\eqref{lambdaT} readily implies that the principal eigenvalue $\lambda_{\delta,m_D,\mu,\rmax}$ is continuous with respect to the parameters $(\delta,m_D,\mu,\rmax)\in(0,+\infty)\times[0,+\infty)\times(0,+\infty)\times\R$.

\vskip 0.3cm
\noindent{\it Step 2: monotonicity with respect to $\delta>0$.} Let us now study the monotonicity and limiting properties of~$\lambda_{\delta,m_D,\mu,\rmax}$ with respect to the parameters $\delta>0$, $m_D\ge0$ and $\mu>0$. Let us start with the dependence with respect to $\delta>0$. To do so, let us fix $m_D\ge0$, $\mu>0$ and~$\rmax\in\R$. or each $\varphi\in\mathcal{E}$, the map $\delta\mapsto\mathcal{R}_{\delta,m_D,\mu,\rmax}(\varphi)$ is nondecreasing in $(0,+\infty)$, from the Cauchy-Schwarz inequality. Therefore, the map $\delta\mapsto\lambda_{\delta,m_D,\mu,\rmax}$ is nondecreasing in~$(0,+\infty)$. \Bk Furthermore, we claim that, when $m_D>0$, the map $\delta\mapsto\lambda_{\delta,m_D,\mu,\rmax}$ is not only nondecreasing but also increasing in $(0,+\infty)$. Indeed, to do so, assume by way of contradiction that there are two migration rates $0<\delta<\delta'$ such that $\lambda_{\delta,m_D,\mu,\rmax}=\lambda_{\delta',m_D,\mu,\rmax}$. Let $\varphi_1^{\delta}$ and $\varphi_1^{\delta'}$ be the functions defined in Lemma~\ref{lem:fpp}, with migration rates $\delta$ and $\delta'$ respectively, and normalized so that~$\|\varphi_1^{\delta}\|_{L^2(\R^n)}=\|\varphi_1^{\delta'}\|_{L^2(\R^n)}=1$. The functions $\varphi_1^{\delta}$ and $\varphi_1^{\delta'}$ respectively minimize $\mathcal{R}_{\delta,\beta,\mu,\rmax}$ and $\mathcal{R}_{\delta',\beta,\mu,\rmax}$ in~$\mathcal{E}$. Thus, the monotonicity of the map~$d\mapsto\mathcal{R}_{d,\beta,\mu,\rmax}(\varphi^{\delta'}_1)$ in $(0,+\infty)$ yields:
$$\lambda_{\delta,m_D,\mu,\rmax}\le\mathcal{R}_{\delta,\beta,\mu,\rmax}(\varphi_1^{\delta'})\le\mathcal{R}_{\delta',\beta,\mu,\rmax}(\varphi_1^{\delta'})=\lambda_{\delta',m_D,\mu,\rmax}=\lambda_{\delta,m_D,\mu,\rmax},$$
hence $\lambda_{\delta,m_D,\mu,\rmax}=\mathcal{R}_{\delta,\beta,\mu,\rmax}(\varphi_1^{\delta'})$, that is,~$\varphi_1^{\delta'}$ also minimizes $\mathcal{R}_{\delta,\beta,\mu,\rmax}$ in $\mathcal{E}$. Therefore, since both functions $\varphi_1^{\delta'}$ and $\varphi_1^{\delta}$ are positive with unit $L^2(\R^n)$ norm, one gets that~$\varphi_1^{\delta'}=\varphi_1^{\delta}$ and:
$$\delta'\,(\varphi_1^{\delta}-\varphi_1^{\delta}\circ\iota)\equiv\delta\,(\varphi_1^{\delta}-\varphi_1^{\delta}\circ\iota)\ \hbox{ in $\R^n$},$$
from the equations satisfied by $\varphi_1^{\delta}=\varphi_1^{\delta'}$. As a consequence, $\varphi_1^{\delta}=\varphi_1^{\delta}\circ\iota$, that is, $\varphi_1^{\delta}=\varphi_2^{\delta}$ by Lemma~\ref{lem:fpp}, where $\varphi_2^\delta$ denotes the function $\varphi_2$ of Lemma~\ref{lem:fpp} (for the migration rate $\delta$). Finally, the system~$\A\,(\varphi_1^{\delta},\varphi_2^{\delta})=\lambda_{\delta,m_D,\mu,\rmax}(\varphi_1^{\delta},\varphi_2^{\delta})$ yields~$r_1\,\varphi_1^{\delta}=r_2\,\varphi_1^{\delta}$ in~$\R^n$, which is clearly impossible since $\varphi_1^{\delta}>0$ in $\R^n$ and $\opt_1\neq\opt_2$ (since one has assumed in the last part of this paragraph that~$m_D>0$). Therefore, the map $\delta\mapsto\lambda_{\delta,m_D,\mu,\rmax}$ is increasing in~$(0,+\infty)$ if $m_D>0$.

Let us now investigate the limits of $\lambda_{\delta,m_D,\mu,\rmax}$ as $\delta\to0$ and $\delta\to+\infty$. First of all, one knows from~\eqref{ineqlambdabis} that $\lambda_{\delta,m_D,\mu,\rmax}\ge-\rmax$ for all $\delta>0$ (this property can also be viewed as a consequence of~\eqref{lambdadelta} since $-r_1(\x)=-\rmax+\|\x-\opt_1\|^2/2\ge-\rmax$ for all $\x\in\R^n$). Furthermore, by choosing a symmetric test function, such as $\varphi_0(x)=\pi^{-n/4}e^{-\|\x\|^2/2}$ for instance, one has $\lambda_{\delta,m_D,\mu,\rmax}\le\mathcal{R}_{\delta,\beta,\mu,\rmax}(\varphi_0)$, and the quantity~$\mathcal{R}_{\delta,\beta,\mu,\rmax}(\varphi_0)$ is independent of $\delta$, hence $\sup_{\delta>0}\lambda_{\delta,m_D,\mu,\rmax}<+\infty$. Therefore, there are two real numbers $\ell^0<\ell^\infty$ in~$[-\rmax,+\infty)$ such that:
$$\lambda_{\delta,m_D,\mu,\rmax}\to\ell^0\hbox{ as $\delta\to0$},\ \hbox{ and }\ \lambda_{\delta,m_D,\mu,\rmax}\to\ell^\infty\hbox{ as $\delta\to+\infty$}.$$

By defining $\mathcal{R}_{0,\beta,\mu,\rmax}(\varphi)$ as above by deleting the (nonnegative) last term of  $\mathcal{R}_{\delta,\beta,\mu,\rmax}(\varphi)$, one has:
$$\mathcal{R}_{0,\beta,\mu,\rmax}(\varphi)\le\mathcal{R}_{\delta,\beta,\mu,\rmax}(\varphi)\le\mathcal{R}_{0,\beta,\mu,\rmax}(\varphi)+2\delta,       
$$
for every function~$\varphi\in\mathcal{E}$. Thus, as $\delta\to0$, the minimum $\lambda_{\delta,m_D,\mu,\rmax}$ of $\mathcal{R}_{\delta,\beta,\mu,\rmax}$ over $\mathcal{E}$ converges to the minimum $\ell^0$ of $\mathcal{R}_{0,\beta,\mu,\rmax}$ over the same set, and this last minimum $\ell^0$ corresponds to the principal eigenvalue of the Schr\"odinger operator,
$$-\frac{\dd^2}{2}\,\Delta-r_1(\x)=-\frac{\dd^2}{2}\,\Delta-\rmax+\frac{\|\x-\opt_1\|^2}{2},$$
acting on the same set of functions. Since the principal eigenvalue of the operator $-\Delta+\|\x\|^2$ is equal to $n$ (with ground state, namely the principal eigenfunction, $\varphi_{GS}(\x)=e^{-\|\x\|^2/2}$ up to multiplicative constants), it easily follows by translation and scaling that:
$$\ell^0=-\rmax+\frac{\dd\,n}{2},$$
with principal eigenfunction~$\varphi(\x)=e^{-\|\x-\opt_1\|^2/(2\mu)}$ up to multiplicative constants.

In order to identity the real number $\ell^\infty=\lim_{\delta\to+\infty}\lambda_{\delta,m_D,\mu,\rmax}=\lim_{k\to+\infty}\lambda_{k,m_D,\mu,\rmax}$, we consider a sequence of (positive) principal eigenfunctions $(\varphi^k_1,\varphi^k_2)_{k\in\N}=(\varphi^k_1,\varphi^k_1\circ\iota)_{k\in\N}$ given by Lemma~\ref{lem:fpp} (with migration rate $\delta=k\in\N$), normalized by $\|\varphi_1^k\|_{L^2(\R^n)}=1$. For each $k\in\N$, there holds $\lambda_{k,m_D,\mu,\rmax}=\mathcal{R}_{k,\beta,\mu,\rmax}(\varphi^k_1)$, hence:
\be\label{ineqk}\baa{r}
\displaystyle\frac{\dd^2}{2}\!\!\int_{\R^n}\!\|\nabla\varphi^k_1(\x)\|^2\,\md\x+\!\int_{\R^n}\!\!\!\frac{\|\x-\opt_1\|^2}{2}\,\varphi^k_1(\x)^2\,\md\x+k\!\!\int_{\R^n}\!\!\big(\varphi^k_1(\x)^2-\varphi^k_1(\x)\varphi^k_1(\iota(\x))\big)\,\md\x\vspace{3pt}\\
=\rmax+\lambda_{k,m_D,\mu,\rmax}.\eaa
\ee
Notice that the right-hand side is bounded as $k\to+\infty$, while the left-hand side is the sum of three nonnegative terms. Therefore, the sequence $(\varphi^k_1)_{k\in\N}$ is bounded in $H^1(\R^n)$ and, up to extraction of a subsequence, there exists a nonnegative function $\varphi_1\in H^1(\R^n)$ such that~$\varphi_1^k\to\varphi_1$ in $L^2_{loc}(\R^n)$ strongly, in~$H^1(\R^n)$ weakly, and almost everywhere in $\R^n$. Furthermore, since~$\|\x-\opt_1\|\to+\infty$ as~$\|\x\|\to+\infty$, one has $\sup_{k\in\N}\|\varphi^k_1\|_{L^2(\R^n\setminus B(0,R))}\to0$ as $R\to+\infty$, hence $\varphi_1^k\to\varphi_1$ in $L^2(\R^n)$ as~$k\to+\infty$, and~$\|\varphi_1\|_{L^2(\R^n)}=1$. Fatou's lemma also implies that the function $\x\mapsto\|\x-\opt_1\|\,\varphi_1(\x)$ belongs to~$L^2(\R^n)$, and so does the function  $\x\mapsto\|\x\|\,\varphi_1(\x)$. Moreover,
$$\int_{\R^n}\big(\varphi^k_1(\x)^2-\varphi^k_1(\x)\varphi^k_1(\iota(\x))\big)\,\md\x\to\int_{\R^n}\big(\varphi_1(\x)^2-\varphi_1(\x)\varphi_1(\iota(\x))\big)\,\md\x,\quad \hbox{ as }k\to+\infty.$$
But since the left-hand side is $O(1/k)$ as $k\to+\infty$ by~\eqref{ineqk}, one gets that:
$$\int_{\R^n}\big(\varphi_1(\x)^2-\varphi_1(\x)\varphi_1(\iota(\x))\big)\,\md\x=0.$$
Since both functions $\varphi_1$ and $\varphi_1\circ\iota$ are nonnegative and with the same $L^2(\R^n)$ norm (equal to~$1$), the case of equality in the Cauchy-Schwarz inequality implies that:
$$\varphi_1=\varphi_1\circ\iota,$$
almost everywhere in $\R^n$. Since each $C^\infty_0(\R^n)$ function $\varphi_1^k+\varphi_2^k=\varphi_1^k+\varphi_1^k\circ\iota$ obeys:
$$-\frac{\dd^2}{2}\,\Delta(\varphi_1^k+\varphi_2^k)-r_1\,\varphi^k_1-r_2\,\varphi^k_2=\lambda_{k,m_D,\mu,\rmax}(\varphi_1^k+\varphi_2^k)\ \hbox{ in }\R^n,$$
and since $\varphi_1^k\to\varphi_1$ and $\varphi_2^k=\varphi_1^k\circ\iota\to\varphi_1\circ\iota=\varphi_1$ in $L^2(\R^n)$ strongly and in $H^1_{loc}(\R^n)$ weakly, it then follows from a passage to the limit in the weak sense and from standard elliptic regularity theory that the function $\varphi_1$ is a $C^\infty(\R^n)$ solution of:
$$-\frac{\dd^2}{2}\,\Delta\varphi_1-\frac{r_1+r_2}{2}\,\varphi_1=\ell^\infty\varphi_1\ \hbox{ in $\R^n$}.$$
Furthermore, since $\|\varphi_1\|_{L^2(\R^n)}=1$ and since $\varphi_1$ is nonnegative, the elliptic strong maximum principle implies that $\varphi_1>0$ in~$\R^n$. The $H^1(\R^n)$ function $\varphi_1$ is then a ground state of the Schr\"odinger operator $-(\dd^2/2)\Delta-(r_1+r_2)/2=-(\dd^2/2)\Delta-\rmax+m_D/4+\|\x\|^2/2$ (remember that $m_D=2\beta^2$). As a consequence,~$\ell^\infty$ is the principal eigenvalue of this operator and $\varphi_1$ is its principal eigenfunction. In other words,
$$\ell^\infty=-\rmax+\frac{\dd\,n}{2}+\frac{m_D}{4},$$
and $\varphi_1(\x)=(\pi\dd)^{-n/4}\,e^{-\|\x\|^2/(2\mu)}$.

\vskip 0.3cm
\noindent{\it Step 3: dependence with respect to $m_D\ge0$.} Let us fix $\delta>0$, $\mu>0$ and $\rmax\in\R$. For any $m_D>0$, one already knows from Step~2 that:
$$\lambda_{\delta,m_D,\mu,\rmax}>-\rmax+\frac{\mu\,n}{2}.$$
By continuity with respect to $m_D\ge0$, one gets that $\lambda_{\delta,0,\mu,\rmax}\ge-\rmax+\mu\,n/2$. Furthermore, the function:
$$\varphi_0:\x\mapsto(\pi\mu)^{-n/4}e^{-\|\x\|^2/(2\mu)},$$
belongs to $\mathcal{E}$ and it is symmetric (that is, $\varphi_0=\varphi_0\circ\iota$ in $\R^n$), hence $\lambda_{\delta,0,\mu,\rmax}\le\mathcal{R}_{\delta,0,\mu,\rmax}(\varphi_0)=-\rmax+\mu\,n/2$. Finally,
$$\lambda_{\delta,0,\mu,\rmax}=-\rmax+\frac{\mu\,n}{2}.$$

Let us now show that $\lambda_{\delta,m_D,\mu,\rmax}\to-\rmax+\mu\,n/2+\delta$ as $m_D\to+\infty$. First of all, for any $m_D=2\beta^2\ge0$, since $\varphi_0(\cdot-\opt_1)\in\mathcal{E}$ (with $\opt_1=(-\beta,0,\dots,0)$ and $\beta\ge0$), one has:
$$\baa{rcl}
\lambda_{\delta,m_D,\mu,\rmax}\le\mathcal{R}_{\delta,\beta,\mu,\rmax}(\varphi_0(\cdot-\opt_1)) & = & \displaystyle-\rmax+\frac{\mu\,n}{2}+\delta-\delta\int_{\R^n}\varphi_0(\x-\opt_1)\,\varphi_0(\iota(\x)-\opt_1)\,\md\x,\vspace{3pt}\\
& < & \displaystyle-\rmax+\frac{\mu\,n}{2}+\delta.\eaa$$
Call $\varphi_1^{m_D}$ the principal eigenfunction given in Lemma~\ref{lem:fpp}. Remember that $\varphi_1^{m_D}$ is positive in~$\R^n$, and let us assume without loss of generality that $\|\varphi_1^{m_D}\|_{L^2(\R^n)}=1$, hence $\varphi_1^{m_D}\in\mathcal{E}$, from the bounds derived in the proof of Lemma~\ref{lem:fpp}. Calling $\psi^{m_D}=\varphi_1^{m_D}(\cdot+\opt_1)\in\mathcal{E}$, one has:
\be\label{psimd}\baa{rcl}
\displaystyle-\rmax+\frac{\mu\,n}{2}+\delta & > & \lambda_{\delta,m_D,\mu,\rmax},\vspace{3pt}\\
& = & \mathcal{R}_{\delta,\beta,\mu,\rmax}(\varphi_1^{m_D}),\vspace{3pt}\\
& = & \displaystyle-\rmax+\frac{\mu^2}{2}\!\!\int_{\R^n}\!\!\|\nabla\varphi_1^{m_D}(\x)\|^2\,\md\x+\!\!\int_{\R^n}\!\!\frac{\|\x\!-\!\opt_1\|^2}{2}\,\varphi_1^{m_D}(\x)^2\,\md\x\vspace{3pt}\\
& & \qquad\qquad \qquad \qquad \displaystyle+\delta\int_{\R^n}\big(\varphi_1^{m_D}(\x)^2-\varphi_1^{m_D}(\x)\,\varphi_1^{m_D}(\iota(\x))\big)\,\md\x,\vspace{3pt}\\
& = & \displaystyle-\rmax+\frac{\mu^2}{2}\!\!\int_{\R^n}\!\!\|\nabla\psi^{m_D}(\x)\|^2\,\md\x+\!\!\int_{\R^n}\!\!\frac{\|\x\|^2}{2}\,\psi^{m_D}(\x)^2\,\md\x\vspace{3pt}\\
& &\qquad\qquad \quad \displaystyle+\delta\int_{\R^n}\big(\psi^{m_D}(\x)^2-\psi^{m_D}(\x)\,\psi^{m_D}(\iota(\x)-2\opt_1)\big)\,\md\x.\eaa
\ee
Since:
$$\min_{\psi\in\mathcal{E}}\Big(\frac{\mu^2}{2}\int_{\R^n}\|\nabla\psi(\x)\|^2\,\md\x+\int_{\R^n}\frac{\|\x\|^2}{2}\,\psi(\x)^2\,\md\x\Big)=\frac{\mu\,n}{2},$$
one obtains that:
\be\label{psimd2}
-\rmax+\frac{\mu\,n}{2}+\delta>\lambda_{\delta,m_D,\mu,\rmax}\ge-\rmax+\frac{\mu\,n}{2}+\delta-\delta\int_{\R^n}\psi^{m_D}(\x)\,\psi^{m_D}(\iota(\x)-2\opt_1)\,\md\x.
\ee
On the other hand, since the last three terms of the last right-hand side of~\eqref{psimd} are nonnegative, one infers from~\eqref{psimd} that:
$$\int_{\R^n}\frac{\|\x\|^2}{2}\,\psi^{m_D}(\x)^2\,\md\x<\frac{\mu\,n}{2}+\delta.$$
Consider any radius $R>0$. The last inequality implies that $\|\psi^{m_D}\|_{L^2(\R^n\setminus B(0,R))}^2<(\mu\,n+2\delta)/R^2$ for every $m_D\ge0$. Now, for every $m_D\ge2R^2$ (that is, $\beta\ge R$), one has $\|\iota(\x)-2\opt_1\|\ge R$ for all $x\in B(0,R)$, hence the Cauchy-Schwarz inequality and the fact that $\|\psi^{m_D}\|_{L^2(\R^n)}=1$ yield:
$$\baa{rcl}
\displaystyle0<\int_{\R^n}\psi^{m_D}(\x)\,\psi^{m_D}(\iota(\x)-2\opt_1)\,\md\x & = &  \displaystyle\int_{B(0,R)}\psi^{m_D}(\x)\,\psi^{m_D}(\iota(\x)-2\opt_1)\,\md\x\vspace{3pt}\\
& &\qquad \qquad \qquad  \displaystyle+\int_{\R^n\setminus B(0,R)}\psi^{m_D}(\x)\,\psi^{m_D}(\iota(\x)-2\opt_1)\,\md\x,\vspace{3pt}\\
& \le & \|\psi^{m_D}\|_{L^2(B(0,R))}\,\|\psi^{m_D}\|_{L^2(\R^n\setminus B(0,R))}\vspace{3pt}\\
& & \quad  +\|\psi^{m_D}\|_{L^2(\R^n\setminus B(0,R))}\,\|\psi^{m_D}(\iota(\cdot)-2\opt_1)\|_{L^2(\R^n\setminus B(0,R))},\vspace{3pt}\\
& \le & \displaystyle 2\,\|\psi^{m_D}\|_{L^2(\R^n\setminus B(0,R))}\le\frac{2\sqrt{\mu\,n+2\delta}}{R}.\eaa$$
Together with~\eqref{psimd2}, it follows that:
$$-\rmax+\frac{\mu\,n}{2}+\delta\ge\limsup_{m_D\to+\infty}\lambda_{\delta,m_D,\mu,\rmax}\ge\liminf_{m_D\to+\infty}\lambda_{\delta,m_D,\mu,\rmax}\ge-\rmax+\frac{\mu\,n}{2}+\delta-\frac{2\delta\sqrt{\mu\,n+2\delta}}{R}.$$
Since $R>0$ can be arbitrarily large, one concludes that:
$$\lambda_{\delta,m_D,\mu,\rmax}\to-\rmax+\frac{\mu\,n}{2}+\delta\ \hbox{ as }m_D\to+\infty.$$

To complete Step~3, let us show that the map $m_D\mapsto\lambda_{\delta,m_D,\mu,\rmax}$ is increasing in $[0,+\infty)$. With the same notations as in the previous paragraph, we claim that, for any $m_D=2\beta^2>0$ (with $\beta>0$):
\be\label{claimpsi}
\left\{\baa{ll}
\varphi_1^{m_D}(x_1,\dots,x_n)\le\varphi_1^{m_D}(-x_1\!-\!2\beta,x_2,\dots,x_n)\vspace{3pt}\\
\varphi_1^{m_D}(x_1,\dots,x_n)\not\equiv\varphi_1^{m_D}(-x_1\!-\!2\beta,x_2,\cdots,x_n)\eaa\right.\hbox{ in }\underbrace{\big\{\x=(x_1,\dots,x_n)\in\R^n,\,x_1\!\le\!-\beta\big\}}_{=:H}.
\ee
Since the proof of this claim is a bit technical, it is postponed below, just before the proof of Theorem~\ref{thm:same_dens_large_delta}. Let us here complete the proof of the monotonicity of the map $m_D\mapsto\lambda_{\delta,m_D,\mu,\rmax}$ in $[0,+\infty)$. Consider any $m_D=2\beta^2>0$, with $\beta>0$. For all $h\in(0,m_D)$, by calling $\beta'=\sqrt{\beta^2-h/2}\in(0,\beta)$, one has:
$$\baa{rcl}
\lambda_{\delta,m_D-h,\mu,\rmax} & \le & \mathcal{R}_{\delta,\beta',\mu,\rmax}(\varphi_1^{m_D}),\vspace{3pt}\\
& = & \displaystyle-\rmax+\frac{\mu^2}{2}\!\!\int_{\R^n}\!\!\|\nabla\varphi_1^{m_D}(\x)\|^2\,\md\x+\!\!\int_{\R^n}\!\!\frac{\|(x_1+\beta',x_2,\dots,x_n)\|^2}{2}\,\varphi_1^{m_D}(\x)^2\,\md\x\vspace{3pt}\\
& & \qquad\qquad\qquad\qquad\qquad\qquad\quad \displaystyle+\delta\int_{\R^n}\big(\varphi_1^{m_D}(\x)^2-\varphi_1^{m_D}(\x)\,\varphi_1^{m_D}(\iota(\x))\big)\,\md\x,\vspace{3pt}\\
& = & \displaystyle\mathcal{R}_{\delta,\beta,\mu,\rmax}(\varphi_1^{m_D})+\frac{(\beta-\beta')^2}{2}-(\beta-\beta')\int_{\R^n}(x_1+\beta)\,\varphi_1^{m_D}(\x)^2\,\md\x,\vspace{3pt}\\
& \le & \displaystyle\lambda_{\delta,m_D,\mu,\rmax}+\frac{(\beta-\beta')^2}{2}\vspace{3pt}\\
& &\qquad \displaystyle-(\beta-\beta')\,\underbrace{\int_{H}(x_1+\beta)\,\big(\varphi_1^{m_D}(\x)^2-\varphi_1^{m_D}(-x_1-2\beta,x_2,\cdots,x_n)^2\big)\,\md\x}_{=:I}.
\eaa$$
From~\eqref{claimpsi} together with the positivity and continuity of $\varphi_1^{m_D}$ in $\R^n$, it follows that the integral~$I$ is positive. Since $I$ does not depend on $h\in(0,m_D)$ and since $\beta'=\sqrt{\beta^2-h/2}$, one infers that:
$$\limsup_{h\to0^+}\frac{\lambda_{\delta,m_D-h,\mu,\rmax}-\lambda_{\delta,m_D,\mu,\rmax}}{h}\le-\frac{I}{4\beta}=-\frac{I}{\sqrt{8m_D}}<0.$$
This above strict inequality is valid for any $m_D>0$. Since the map $m_D\mapsto\lambda_{\delta,m_D,\mu,\rmax}$ is continuous in $[0,+\infty)$, one then concludes that it is increasing in $[0,+\infty)$.

\vskip 0.3cm
\noindent{\it Step 4: monotonicity with respect to $\mu>0$.} Let us fix here $\delta>0$, $m_D=2\beta^2\ge0$  and $\rmax\in\R$. Remember that $\lambda_{\delta,m_D,\mu,\rmax}=\min_{\varphi\in\mathcal{E}}\mathcal{R}_{\delta,\beta,\mu,\rmax}(\varphi)$ and that the minimum is reached only by~$\pm\varphi_1$, where $\varphi_1$ is the principal eigenfunction given in Lemma~\ref{lem:fpp}, normalized with unit $L^2(\R^n)$ norm. But $\mathcal{R}_{\delta,\beta,\mu,\rmax}(\varphi)$ is nondecreasing with respect to $\mu>0$ for each $\varphi\in\mathcal{E}$, and the principal eigenfunctions $\varphi_1$ are non constant (that is, the $L^2(\R^n)$ norm of their gradient is positive). One then infers that $\lambda_{\delta,m_D,\mu,\rmax}$ is increasing with respect to $\mu$.

Notice now that $\lambda_{\delta,m_D,\mu,\rmax}=\lambda_{\delta,m_D,\mu,0}-\rmax\ge-\rmax$, and call:
$$\lambda^0=\lim_{\mu\to0^+}\lambda_{\delta,m_D,\mu,0}\ge0,$$
that is, $\lambda_{\delta,m_D,\mu,\rmax}\to-\rmax+\lambda^0$ as $\mu\to0^+$.

Let us show in this paragraph that $\lambda^0\le\min(\delta,m_D/4)$. First of all, consider a $C^1(\R^n)$ radially symmetric function $\varphi$ with compact support and unit $L^2(\R^n)$ norm. For $\varepsilon>0$ and $\x\in\R^n$, call $\varphi_\varepsilon(\x)=\epsilon^{-n/2}\varphi(\x/\varepsilon)$. Each function $\varphi_\varepsilon$ is radially symmetric and belongs to $\mathcal{E}$, hence:
$$\lambda_{\delta,m_D,\mu,0}\le\mathcal{R}_{\delta,\beta,\mu,0}(\varphi_\varepsilon)=\frac{\mu^2}{2}\int_{\R^n}\|\nabla\varphi_\varepsilon(\x)\|^2\,\md\x+\int_{\R^n}\frac{\|\x-\opt_1\|^2}{2}\,\varphi_\varepsilon(\x)^2\,\md\x,$$
and:
$$\lambda^0\le\int_{\R^n}\frac{\|\x-\opt_1\|^2}{2}\,\varphi_\varepsilon(\x)^2\,\md\x,$$
at the limit $\mu\to0^+$. Since the above inequality holds for all $\varepsilon>0$ and since the right-hand side converges to $\|\opt_1\|^2/2=m_D/4$ as $\varepsilon\to0^+$, one gets that:
$$\lambda^0\le\frac{m_D}{4}.$$
When $m_D=0$, then $\lambda^0=0$. Assume in the sequel that $m_D>0$. Each function $\psi_\varepsilon:=\varphi_\varepsilon(\cdot-\opt_1)$ belongs to $\mathcal{E}$ as well, hence:
$$\baa{rcl}
\lambda_{\delta,m_D,\mu,0} & \le & \mathcal{R}_{\delta,\beta,\mu,0}(\psi_\varepsilon),\vspace{3pt}\\
& = & \displaystyle\frac{\mu^2}{2}\int_{\R^n}\|\nabla\psi_\varepsilon(\x)\|^2\,\md\x+\int_{\R^n}\frac{\|\x-\opt_1\|^2}{2}\,\psi_\varepsilon(\x)^2\,\md\x+\delta-\delta\int_{\R^n}\psi_\varepsilon(\x)\,\psi_\varepsilon(\iota(\x))\,\md\x,\vspace{3pt}\\
& = & \displaystyle\frac{\mu^2}{2}\int_{\R^n}\|\nabla\varphi_\varepsilon(\x)\|^2\,\md\x+\int_{\R^n}\frac{\|\x\|^2}{2}\,\varphi_\varepsilon(\x)^2\,\md\x+\delta-\delta\int_{\R^n}\varphi_\varepsilon(\x)\,\varphi_\varepsilon(\iota(\x)-2\opt_1)\,\md\x,\eaa$$
and:
$$\lambda^0\le\int_{\R^n}\frac{\|\x\|^2}{2}\,\varphi_\varepsilon(\x)^2\,\md\x+\delta-\delta\int_{\R^n}\varphi_\varepsilon(\x)\,\varphi_\varepsilon(\iota(\x)-2\opt_1)\,\md\x,$$
at the limit $\mu\to0^+$. Since the above inequality holds for all $\varepsilon>0$ and since the right-hand side converges to $0+\delta-0=\delta$ as $\varepsilon\to0^+$ (since $\varphi$ has compact support, and $\|\opt_1\|>0$), one gets that:
$$\lambda^0\le\delta.$$
To sum up, $\lambda_0\in[0,\min(\delta,m_D/4)]$ for all $\delta>0$ and $m_D\ge0$. 

It only remains to prove that $\lambda_{\delta,m_D,\mu,\rmax}\to+\infty$ as $\mu\to+\infty$. Since the map $\mu\mapsto\lambda_{\delta,m_D,\mu,\rmax}$ is increasing and since $\lambda_{\delta,m_D,\mu,\rmax}\ge-\rmax$ for all $\mu>0$, there is $\lambda^\infty\in(-\rmax,+\infty]$ such that $\lambda_{\delta,m_D,\mu,\rmax}\displaystyle\mathop{\to}^<\lambda^\infty$ as $\mu\to+\infty$. Let $\varphi_1^k$ be the principal eigenfunction given in Lemma~\ref{lem:fpp} with mutational parameter $\mu=k\in\N$ with $k\ge1$, and normalized with unit $L^2(\R^n)$ norm. Hence,~$\varphi_1^k\in\mathcal{E}$, and:
$$\baa{r}
\displaystyle\frac{k^2}{2}\int_{\R^n}\|\nabla\varphi_1^k(\x)\|^2\,\md\x-\rmax+\int_{\R^n}\frac{\|\x-\opt_1\|^2}{2}\,\varphi_1^k(\x)^2\,\md\x+\delta\int_{\R^n}\big(\varphi_1^k(\x)^2-\varphi_1^k(\x)\,\varphi_1^k(\iota(\x))\big)\,\md\x\vspace{3pt}\\
=\mathcal{R}_{\delta,\beta,k,\rmax}(\varphi_1^k)=\lambda_{\delta,m_D,k,\rmax}<\lambda^\infty,\eaa$$
for all $k\in\N$ with $k\ge1$. Since the last term of the left-hand side is nonnegative by the Cauchy-Schwarz inequality, one infers that:
\be\label{psikmu}
\frac{k^2}{2}\int_{\R^n}\|\nabla\psi^k(\x)\|^2\,\md\x+\int_{\R^n}\frac{\|\x\|^2}{2}\,\psi^k(\x)^2\,\md\x<\rmax+\lambda^\infty,
\ee
for all $k\in\N$ with $k\ge1$, with $\psi^k:=\varphi_1^k(\cdot+\opt_1)\in\mathcal{E}$. Assume now by way of contradiction that~$\lambda^\infty<+\infty$ and choose $R>0$ such that $2(\rmax+\lambda^\infty)/R^2\le3/4$. Hence,
$$\frac{R^2}{2}\int_{\R^n\setminus B(0,R)}\psi^k(\x)^2\,\md\x<\rmax+\lambda^\infty\le\frac34\times\frac{R^2}{2},$$
and $1\!\ge\!\|\psi^k\|_{L^2(B(0,R))}\!=\!\sqrt{\int_{B(0,R)}\!\psi^k(\x)^2\md\x}\ge\!1/2$ for all $k\!\in\!\N$ with $k\ge1$ (remember that~$\|\psi^k\|_{L^2(\R^n)}\!=\!1$). The inequality~\eqref{psikmu} also implies that the sequence $(\psi^k)_{k\in\N^*}$ is bounded in $H^1(\R^n)$ and that~$\|\nabla\psi^k\|_{L^2(\R^n)}\to0$ as~$k\to+\infty$. Up to extraction of a subsequence, there is a function~$\psi^\infty\in H^1(\R^n)$ such that $\psi^k\to\psi^\infty$ in~$H^1(\R^n)$ weakly, and in $L^2_{loc}(\R^n)$ strongly. In particular, $1/2\le\|\psi^\infty\|_{L^2(B(0,R))}\le1$. Furthermore, $\|\nabla\psi^\infty\|_{L^2(\R^n)}=0$, that is, there is a cons\-tant~$C$ such that~$\psi^\infty=C$ almost everywhere in $\R^n$. But this constant $C$ can not be zero since~$\|\psi^\infty\|_{L^2(B(0,R))}>0$, and then $\psi^\infty$ can not be in $H^1(\R^n)$. One has then reached a contradiction, hence $\lambda^\infty=+\infty$ and:
\[\lambda_{\delta,m_D,\mu,\rmax}\to+\infty\ \hbox{ as $\mu\to+\infty$}.\]
The proof of Proposition~\ref{prop:continuity_lambda} is thereby complete.
\qed

\vskip0.3cm

\noindent \textit{Proof of Eq.~\eqref{claimpsi}.} Throughout this proof, we fix $(\delta,\mu,\rmax)\in(0,+\infty)\times(0,+\infty)\times\R$, as well as $m_D=2\beta^2>0$ with $\beta>0$ and $\opt_1=(-\beta,0,\dots,0)$.

Let us first show that the function $\varphi_1^{m_D}$ and the reflected one $\x\mapsto\varphi_1^{m_D}(-x_1-2\beta,x_2,\dots,x_n)$ can not be identically equal in $H$. If they were, then they would be identically equal in $\R^n$ by definition of $H$, where $H$ is the half-space defined in~\eqref{claimpsi}. From the equations satisfied by these two functions, it easily follows that $\varphi_1^{m_D}(\x)= \varphi_1^{m_D}(\x+4\opt_1)$ for all $\x\in\R^n$. In other words,~$\varphi_1^{m_D}$ would be periodic, which is ruled out since $\varphi_1^{m_D}$ is a non-trivial function in $H^1(\R^n)$. Therefore,~$\varphi_1^{m_D}$ and $\x\mapsto\varphi_1^{m_D}(-x_1-2\beta,x_2,\dots,x_n)$ can not be identically equal in $H$

It then remains to show the inequality in~\eqref{claimpsi}. To do so, from the proof of Lemma~\ref{lem:fpp}, it is sufficient to show that, for any $R>\beta$, one has:
\be\label{claim1}
\phi_1(\x)\le\phi_1(-x_1-2\beta,x_2,\dots,x_n)\hbox{ for all $\x\in\overline{B(0,R)}$ with $x_1\le-\beta$},
\ee
where $\phi_1$ here denotes the first component of the pair of principal eigenfunctions $(\phi_1,\phi_2)$ of the operator $\mathcal{A}$ defined in~\eqref{eq:defA} with Dirichlet boundary conditions on $\partial B(0,R)$ (with principal eigenvalue denoted $\lambda^R$). Remember that the $C^\infty_0(\overline{B(0,R)})$ functions $\phi_1$ and $\phi_2$ are positive in~$B(0,R)$ and solve $\mathcal{A}(\phi_1,\phi_2)-\lambda^R(\phi_1,\phi_2)=(0,0)$ in $\overline{B(0,R)}$.

Let us then fix $R>\beta$ till the end of the proof. To show the inequality~\eqref{claim1}, we will actually prove the following stronger property:
\be\label{claim2}
\forall\,\nu\in(-R,-\beta],\ \forall\,i\in\{1,2\},\ \phi_i\le\phi_i^\nu\hbox{ in }\overline{H_\nu},
\ee
where:
$$H_\nu=\{\x=(x_1,\dots,x_n)\in B(0,R),x_1<\nu\},$$
and:
$$\phi_i^\nu(\x)=\phi_i(-x_1+2\nu,x_2,\dots,x_n),$$
for $i=1,2$ (the desired inequality~\eqref{claim1} then follows from~\eqref{claim2} with $i=1$ and $\nu=-\beta$). The proof of~\eqref{claim2} is based on the method of moving planes~\cite{Ale62,GidNiNir79}: it is first proven for $\nu$ larger but close to~$-R$, and then up to the value $-\beta$ by increasing $\nu$ from $-R$ to $-\beta$. Two main ingredients will be used in the proof of~\eqref{claim2}. One of these ingredients is the strong maximum principle (see e.g.~\cite[Proposition~12.1]{BusSir04}) applied to the operator $\mathcal{A}-\lambda^R$: it says that if $\omega$ is an open connected set of $\R^n$ and if $(\psi_1,\psi_2)$ is a pair of nonnegative $C^2(\omega)$ functions solving $\mathcal{A}(\psi_1,\psi_2)-\lambda^R(\psi_1,\psi_2)\ge(0,0)$ componentwise in $\omega$, then either $(\psi_1,\psi_2)\equiv(0,0)$ in $\omega$, or both functions $\psi_1$ and $\psi_2$ are positive in $\omega$. The second main ingredient is the weak maximum principle in subsets of $B(0,R)$ with small Lebesgue measure~\cite[Corollary~14.1]{BusSir04}: it says that there is $\eta>0$ such that, if $\omega$ is an open subset of $B(0,R)$ with Lebesgue measure less than $\eta$ and if $(\psi_1,\psi_2)$ is a pair of $C^2(\omega)\cap C(\overline{\omega})$ functions solving $\mathcal{A}(\psi_1,\psi_2)-\lambda^R(\psi_1,\psi_2)\ge(0,0)$ componentwise in $\omega$ and $(\psi_1,\psi_2)\ge(0,0)$ on~$\partial\omega$, then~$(\psi_1,\psi_2)\ge(0,0)$ in $\overline{\omega}$. Before putting these ingredients together, let us first observe that, for every $\nu\in(-R,-\beta]$ and for every $\x\in\overline{H_\nu}$, there holds:
$$|x_1+\beta|\ge|\!-x_1+2\nu+\beta|\ \hbox{ and }\ |x_1-\beta|\ge|\!-x_1+2\nu-\beta|,$$
hence:
$$\mathcal{A}(\phi_1^\nu,\phi_2^\nu)-\lambda^R(\phi_1^\nu,\phi_2^\nu)\ge(0,0)\ \hbox{ in }\overline{H_\nu},$$
from the definitions~\eqref{def:ri} and~\eqref{def:opt_beta} of the fitnesses $r_i$ and the optima $\opt_i$. Therefore, for every $\nu\in(-R,-\beta]$, the $C^2(\overline{H_\nu})$ functions $\phi_1^\nu-\phi_1$ and $\phi_2^\nu-\phi_2$ satisfy:
\be\label{ineqphinu}
\mathcal{A}(\phi_1^\nu-\phi_1,\phi_2^\nu-\phi_2)-\lambda^R(\phi_1^\nu-\phi_1,\phi_2^\nu-\phi_2)\ge(0,0)\ \hbox{ in }\overline{H_\nu}.
\ee
Moreover, for every $\nu\in(-R,-\beta]$, both functions $\phi_1^\nu-\phi_1$ and $\phi_2^\nu-\phi_2$ are nonnegative and not identically equal to $0$ on $\partial H_\nu$ since $-R<\nu\le-\beta<0$ and since $\phi_1$ and $\phi_2$ are positive in~$B(0,R)$ and vanish on $\partial B(0,R)$ (in particular, the functions $\phi_1^\nu-\phi_1$ and $\phi_2^\nu-\phi_2$ can not be identically equal to~$0$ in~$\overline{H_\nu}$). Since the Lebesgue measure of $H_\nu$ goes to $0$ as $\nu\to-R$, the aforementioned weak maximum principle in subsets of $B(0,R)$ with small Lebesgue measure yields the existence of $\nu_0\in(-R,-\beta)$ such that~\eqref{claim2} holds for all $\nu\in(-R,\nu_0]$. Denote now:
$$\nu^*=\sup\big\{\nu\in(-R,-\beta],\,(\phi_1,\phi_2)\le(\phi_1^{\nu'},\phi_2^{\nu'})\hbox{ in }\overline{H_{\nu'}}\hbox{ for all }\nu'\in(-R,\nu]\big\}.$$
One has $-R<\nu_0\le\nu^*\le-\beta$, and one claims that $\nu^*=-\beta$. Assume not. Then $\nu^*<-\beta$. By continuity, one has $(\phi_1,\phi_2)\le(\phi_1^{\nu^*},\phi_2^{\nu^*})$ in $\overline{H_{\nu^*}}$. Remember also that $\phi_1^{\nu^*}-\phi_1$ and $\phi_2^{\nu^*}-\phi_2$ can not be identically equal to~$0$ in~$\overline{H_{\nu^*}}$, and then can not be identically equal to~$0$ in~$H_{\nu^*}$ by continuity. Together with~\eqref{ineqphinu}, the aforementioned strong maximum principle implies that:
$$\phi_1^{\nu^*}-\phi_1>0\hbox{ and }\phi_2^{\nu^*}-\phi_2>0\hbox{ in $H_{\nu^*}$}.$$
Pick a compact subset $K$ of $H_{\nu^*}$ such that the Lebesgue measure of $H_{\nu^*}\setminus K$ is less than $\eta/2$ (where $\eta>0$ is given above, for which the weak maximum principle holds for $\mathcal{A}-\lambda^R$ in open subsets of $B(0,R)$ of measure less than $\eta$). By continuity, $\min_K(\phi_1^{\nu^*}-\phi_1)>0$ and $\min_K(\phi_2^{\nu^*}-\phi_2)>0$, and there is $\varepsilon\in(0,-\beta-\nu^*)$ such that:
$$\min_K(\phi_1^{\nu}-\phi_1)>0\hbox{ and }\min_K(\phi_2^{\nu}-\phi_2)>0\hbox{ for all $\nu\in[\nu^*,\nu^*+\varepsilon]$}.$$
Without loss of generality, one can also assume that the Lebesgue measure of $H_{\nu^*+\varepsilon}\setminus H_{\nu^*}$ is less than $\eta/2$, hence the Lebesgue measure of $H_{\nu}\!\setminus\!K$ is less than $\eta$ for all $\nu\in[\nu^*,\nu^*+\varepsilon]$. Furthermore, $(\phi_1^{\nu}-\phi_1,\phi_2^{\nu}-\phi_2)\ge(0,0)$ on $\partial(H_{\nu}\!\setminus\!K)$ for all $\nu\in[\nu^*,\nu^*+\varepsilon]$, and together with~\eqref{ineqphinu} the aforementioned weak maximum principle then implies that $(\phi_1^{\nu}-\phi_1,\phi_2^{\nu}-\phi_2)\ge(0,0)$ in~$\overline{H_{\nu}\setminus K}$, for all $\nu\in[\nu^*,\nu^*+\varepsilon]$. Finally, $(\phi_1^{\nu}-\phi_1,\phi_2^{\nu}-\phi_2)\ge(0,0)$ in $\overline{H_{\nu}}$ for all $\nu\in[\nu^*,\nu^*+\varepsilon]$, contradicting the definition of $\nu^*$. As a consequence, $\nu^*=-\beta$ and~\eqref{claim2} has been proven for all $\nu\in(-R,-\beta)$ and then also for $\nu=-\beta$ by continuity. As already emphasized, this yields~\eqref{claim1} and then~\eqref{claimpsi}.
\qed

\vskip0.3cm

\noindent \textit{Proof of Theorem~\ref{thm:same_dens_large_delta}}.
Let $\densbf_\delta=(\dens_{\delta,1},\dens_{\delta,2})$ be the unique $C^{1,2}(\R_+\times\R^n)^2$ solution of~\eqref{eq:sys_main1} given by Theorem~$\ref{thm:existunidens}$ and Remark~$\ref{rem:nonsym}$, for growth functions $f_1,f_2$ of the first type~\eqref{eq:f_type1}, with a fixed initial condition $\densbf^0=(u^0_1,u^0_2)$ independent of $\delta$ and such that both functions $u^0_1,u^0_2$ satisfy the assumptions~{\rm{(H1)-(H3)}}. Let us fix two positive times $0<T'\le T$ and let us show that $\sup_{t\in[T',T]}\|u_{\delta,1}(t,\cdot)-u_{\delta,2}(t,\cdot)\|_{L^\infty(\R^n)}\to0$ as $\delta\to+\infty$.

From the first part of the proof of Theorem~\ref{thm:existunidens}, especially from~\eqref{eq:maj1:dens},~\eqref{eq:ineghi}-\eqref{defzeta} and similar calculations as the ones between~\eqref{eq:ineghi} and~\eqref{defzeta}, it follows that there exists a constant $K\ge 0$ (independent of~$\delta>0$) such that, for all $\delta>0$,
\be\label{ineqKdelta}
|x_1\dens_{\delta,2}(t,\x)| \le |x_1h_2(t,\x)|\le K, \quad \hbox{ for all $t\in[0,T]$ and $\x=(x_1\tpp x_n)\in \R^ n$},
\ee
with $h_2$ defined by~\eqref{dfn:hbf} (notice that the function $h_2$ actually depends on $\delta$, but the upper bound~\eqref{eq:ineghi} is independent of $\delta>0$). For each $\delta>0$, one infers from~\eqref{eq:sys_main1}-\eqref{def:ri} and~\eqref{def:opt_beta} that the function $v_\delta:=u_{\delta,1}-u_{\delta,2}$ is a classical $C^{1,2}(\R_+\times\R^n)$ solution of:
\[
\partial_t v_\delta(t,\x)= {\dd^2\over 2} \Delta v_\delta (t,\x) +r_1(\x)\,v_\delta (t,\x)-2\delta\,v_\delta(t,\x) -2\beta x_1\dens_{\delta,2}(t,\x),
\]
such that $v_\delta$ is locally bounded in time and $v_\delta(t,x)\to0$ as $\|\x\|\to+\infty$ locally uniformly in $t\in\R_+$. The previous relation, together with~\eqref{def:m(x)} and~\eqref{ineqKdelta}, implies that:\small
\[
\left\{ \begin{array}{l}
    \displaystyle-2\beta K \le \partial_t v_\delta(t,\x)- \frac{\dd^2}{2} \Delta v_\delta (t,\x) -\big(\rmax + m_1(\x)\big)\, v_\delta (t,\x)+2\delta\,v_\delta (t,\x) \le 2\beta K, \ \ t\in[0,T],\ \x\in \R^n,\vspace{3pt}
  \\
    |v_\delta(0,\x)|\le\max(\|u^0_1\|_{L^\infty(\R^n)},\|u^0_2\|_{L^\infty(\R^n)})=:M, \ \ \x \in \R^n.
\end{array}\right.
\]
\normalsize

Since the potential $m_1(\x)=-\|\x-\opt_1\|^2/2$ is nonpositive, there exists a $C^{1,2}(\R_+\times\R^n)$ solution~$V:\R_+\times\R^n\to[0,M]$ of:
\[
\left\{ \begin{array}{l}
\displaystyle\partial_t V(t,\x)= {\dd^2\over 2} \Delta V(t,\x) +m_1(\x)\,V(t,\x), \ \ t\ge0,\ \x\in \R^n,\vspace{3pt}
  \\
    V(0,\x)=M, \ \ \x \in \R^n.
\end{array}\right.
\]
Such a function $V$, which is independent of $\delta>0$, can be obtained as the nondecreasing local limit as $R\to+\infty$ of $C^{1,2}(\R_+\times\overline{B(0,R)})$ solutions $V^R:\R_+\times\overline{B(0,R)}\to[0,M]$ of the same equation in $\R_+\times\overline{B(0,R)}$, with Dirichlet boundary conditions $V^R=0$ on~$\R_+\times\partial B(0,R)$ and initial conditions of the type~$V^R(0,\x)=M\,\phi(\|\x\|/R)$ in $\overline{B(0,R)}$, where $\phi:[0,1]\to[0,1]$ is a~$C^\infty([0,1])$ nonincreasing function such that $\phi=1$ in $[0,1/3]$ and $\phi=0$ in $[2/3,1]$.

Consider now any $\delta>\rmax/2$ and let $V_\delta$ be the $C^{1,2}(\R_+\times\R^n)$ function defined in $\R_+\times\R^n$ by:
$$V_\delta(t,\x)=v_\delta(t,\x)\,e^{(2\delta-\rmax)t}-\frac{2\beta K}{2\delta-\rmax}\,\big(e^{(2\delta-\rmax)t}-1\big).$$
A straightforward calculation shows that:
$$\partial_tV_\delta(t,\x)-\frac{\dd^2}{2}\,\Delta V_\delta(t,\x)-m_1(\x)\,V_\delta(t,\x)\le m_1(\x)\,\frac{2\beta K}{2\delta-\rmax}\,\big(e^{(2\delta-\rmax)t}-1\big)\le0,$$
for all $(t,x)\in[0,T]\times\R^n$. Furthermore, $V_\delta(0,\x)=v_\delta(0,\x)\le M=V(0,\x)$ for all $\x\in\R^n$, and~$\limsup_{\|\x\|\to+\infty}V_\delta(t,\x)\le0$ uniformly in $t\in[0,T]$. It follows from the maximum principle that $V_\delta(t,\x)\le V(t,x)$ for all $(t,\x)\in[0,T]\times\R^n$, hence:
$$v_\delta(t,\x)\le e^{(\rmax-2\delta)t}V(t,\x)+\frac{2\beta K}{2\delta-\rmax}\,\big(1-e^{(\rmax-2\delta)t}\big),\quad  \hbox{ for all $(t,\x)\in[0,T]\times\R^n$}.$$
Since the function $V$ is bounded (by $M$), one gets that:
$$\limsup_{\delta\to+\infty}\Big(\sup_{[T',T]\times\R^n}v_\delta\Big)\le0,$$
recalling that $0<T'\le T$. The same argument applied to the functions $-V_\delta$ and $-v_\delta$ implies that, for all~$\delta>\rmax/2$ and $(t,\x)\in[0,T]\times\R^n$,
$$v_\delta(t,\x)\ge-e^{(\rmax-2\delta)t}V(t,\x)-\frac{2\beta K}{2\delta-\rmax}\,\big(1-e^{(\rmax-2\delta)t}\big),$$
hence $\liminf_{\delta\to+\infty}\big(\inf_{[T',T]\times\R^n}v_\delta\big)\ge0$. As a conclusion, $\sup_{[T',T]\times\R^n}|v_\delta|\to0$ as $\delta\to+\infty$ and the proof of Theorem~\ref{thm:same_dens_large_delta} is thereby complete.
\qed

\vskip 0.5cm
\noindent{\bf{Acknowledgements.}} The authors are grateful to the anonymous referees for their valuable comments and suggestions, which led to significant improvements of the manuscript.

\bibliographystyle{plain}
\bibliography{biblio2.bib}

\end{document}